\documentclass[12pt,oneside,reqno]{amsart}
\usepackage{mathrsfs}
\usepackage{graphics}
\usepackage{amssymb}
\usepackage{enumerate}
\newcommand{\dif}{\mathrm{d}}

\newcommand{\be}{\begin{eqnarray}}
\newcommand{\ee}{\end{eqnarray}}
\newcommand{\ce}{\begin{eqnarray*}}
\newcommand{\de}{\end{eqnarray*}}
\newtheorem{theorem}{Theorem}[section]
\newtheorem{lemma}[theorem]{Lemma}
\newtheorem{remark}[theorem]{Remark}
\newtheorem{definition}[theorem]{Definition}
\newtheorem{proposition}[theorem]{Proposition}
\newtheorem{Example}[theorem]{Example}
\newtheorem{corollary}[theorem]{Corollary}
\def\e{\varepsilon}
\def\s{\sigma}
\def\t{\theta}
\def\a{\alpha}

\def\b{\beta}
\def\d{\delta}

\def\g{\gamma}
\def\l{\lambda}

\def\[{{\Big[}}
\def\]{{\Big]}}
\def\<{{\langle}}
\def\>{{\rangle}}
\def\({{\Big(}}
\def\){{\Big)}}

\def\no{\nonumber}
\def\bt{\begin{theorem}}
\def\et{\end{theorem}}
\def\bl{\begin{lemma}}
\def\el{\end{lemma}}
\def\br{\begin{remark}}
\def\er{\end{remark}}
\def\bx{\begin{Example}}
\def\ex{\end{Example}}
\def\bd{\begin{definition}}
\def\ed{\end{definition}}
\def\bp{\begin{proposition}}
\def\ep{\end{proposition}}
\def\bc{\begin{corollary}}
\def\ec{\end{corollary}}

\def\cB{{\mathcal B}}

\def\cN{{\mathcal N}}

\def\cP{{\mathcal P}}

\def\mE{{\mathbb E}}

\def\mP{{\mathbb P}}

\def\mR{{\mathbb R}}

\def\mW{{\mathbb W}}

\def\sB{{\mathscr B}}

\def\sF{{\mathscr F}}

\def\sL{{\mathscr L}}

\def\geq{\geqslant}
\def\leq{\leqslant}

\begin{document}

\allowdisplaybreaks

\title{Strong approximation of nonlinear filtering for multiscale McKean-Vlasov stochastic systems}

\author{Huijie Qiao and Wanlin Wei }

\dedicatory{School of Mathematics,
Southeast University,\\
Nanjing, Jiangsu 211189, P.R.China}



\thanks{This work was partly supported by NSF of China (No.12071071).}

\thanks{Corresponding author: Huijie Qiao, hjqiaogean@seu.edu.cn}


\date{}

\begin{abstract}
This work concerns the nonlinear filtering problem of multiscale McKean-Vlasov stochastic systems where the whole systems depend on distributions of fast components. First of all, we prove that the slow component of the original system converges to an average system in the $L^{2p}$ ($p\geqslant 1$) sense. Moreover, we obtain the strong convergence order for the $L^2$ case. Then, given an observation process which depends on the slow component and its distribution, we show that the nonlinear filtering of the slow component and its distribution also converges to that of the average system in the $L^{q}$  ($p\geq 8, 1\leq q\leq \frac{p}{8}$) sense.
\end{abstract}

\maketitle \rm

\section{Introduction}
McKean-Vlasov stochastic differential equations (SDEs for short) are also called distri\\bution-dependent SDEs, or mean-field SDEs. And the difference between McKean-Vlasov SDEs and classical SDEs is that the former depends on the positions and probability distributions of particles. Therefore, McKean-Vlasov SDEs can better describe many models. The study on McKean-Vlasov SDEs was initiated by H. P. McKean \cite{Mckean} who was inspired by Kac's program in Kinetic Theory. Nowadays, McKean-Vlasov SDEs have been widely applied in many fields, such as biology, game theory, optimal control theory and interacting particle systems. And there are many results about them. Let us recall some works. Sznitman proved the existence and uniqueness of strong solutions to McKean-Vlasov SDEs under global Lipschitz conditions in \cite{Sznit}. Ding and Qiao \cite{DQ,DQ2} studied the well-posedness and stability of solutions to McKean-Vlasov SDEs with non-Lipschitz coefficients. Wang investigated the exponential ergodicity of the strong solutions to Landau type McKean-Vlasov SDEs in \cite{WangFY}. Sen and Caines \cite{NPE} and Liu and Qiao \cite{MQ} studied nonlinear filtering problems of McKean-Vlasov SDEs with independent noises and correlated noises, respectively.

Besides, multiscale SDEs, or slow-fast systems, are widely used in engineering and science fields (c.f. \cite{{Goggin1}, {Goggin2}}). The average principle for them was first studied by Khasminskii \cite{rzk}, see \cite{LRSX, LD1, SunX, XLLM} (and the references therein) for further generalizations. Here we briefly mention some results related with ours. Liu \cite{LD1} studied SDEs with two well-separated time scales under Lipschitz conditions and established that the slow part of the original system converges to an average system in the $L^{2}$ sense. Liu, R\"{o}ckner, Sun and Xie considered a class of slow-fast SDEs and proved the convergence in the $L^{p}$ sense in \cite{LRSX}.

Now, multiscale McKean-Vlasov SDEs have also been studied. For example, R\"{o}ckner, Sun and Xie \cite{SunX} investigated the following multiscale McKean-Vlasov SDEs: for $T>0$
\be\left\{\begin{array}{l}
\mathrm{d} \tilde{X}_t^{\varepsilon}=\tilde{b}_1\left(\tilde{X}_t^{\varepsilon}, \mathscr{L}_{\tilde{X}_t^{\varepsilon}}^{\mathbb{P}}, \tilde{Z}_t^{\varepsilon}\right) \mathrm{d} t+\tilde{\sigma}_1\left(\tilde{X}_t^{\varepsilon}, \mathscr{L}_{\tilde{X}_t^{\varepsilon}}^{\mathbb{P}}\right) \mathrm{d} B_t, \\
\tilde{X}_0^{\varepsilon}=x_0, \quad 0 \leqslant t \leqslant T \\
\mathrm{~d} \tilde{Z}_t^{\varepsilon}=\frac{1}{\varepsilon} \tilde{b}_2\left(\tilde{X}_t^{\varepsilon}, \mathscr{L}_{\tilde{X}_t^{\varepsilon}}^{\mathbb{P}}, \tilde{Z}_t^{\varepsilon}\right) \mathrm{d} t+\frac{1}{\sqrt{\varepsilon}} \tilde{\sigma}_2\left(\tilde{X}_t^{\varepsilon}, \mathscr{L}_{\tilde{X}_t^{\varepsilon}}^{\mathbb{P}}, \tilde{Z}_t^{\varepsilon}\right) \mathrm{d} W_t, \\
\tilde{Z}_0^{\varepsilon}=z_0, \quad 0 \leqslant t \leqslant T,
\end{array}\right.
\label{Eq0}
\ee
where $\left(B_t\right),\left(W_t\right)$ are $n$- and $m$-dimensional standard Brownian motions defined on the complete filtered probability space $(\Omega, \mathscr{F},\left\{\mathscr{F}_t\right\}_{t \in[0, T]}, \mathbb{P})$, respectively, these mappings $\tilde{b}_1$ : $\mathbb{R}^n \times \mathcal{P}_2\left(\mathbb{R}^n\right) \times \mathbb{R}^m \rightarrow \mathbb{R}^n, \tilde{\sigma}_1: \mathbb{R}^n \times \mathcal{P}_2\left(\mathbb{R}^n\right) \rightarrow \mathbb{R}^{n \times n}, \tilde{b}_2: \mathbb{R}^n \times \mathcal{P}_2\left(\mathbb{R}^n\right) \times \mathbb{R}^m \rightarrow \mathbb{R}^m$, $\tilde{\sigma}_2: \mathbb{R}^n \times \mathcal{P}_2\left(\mathbb{R}^n\right) \times \mathbb{R}^m \rightarrow \mathbb{R}^{m \times m}$ are all Borel measurable, and $\mathscr{L}_{X_t^{\varepsilon}}^{\mathbb{P}}$ denotes the distribution of $X_t^{\varepsilon}$ under $\mathbb{P}$. There they showed that the slow part converges to an average system in the $L^2$ sense. Later, Hong, Li and Liu \cite{hll} generalized this result to the infinite dimensional space. Very recently, the first named author \cite{Qiao0} also obtain the same result for multiscale multivalued McKean-Vlasov SDEs.

Note that the whole system (\ref{Eq0}) doesn't depend on the distribution of the fast component. One improvement is that Gao, Hong and Liu \cite{ghl} added the distribution of the fast component to the slow component in the system (\ref{Eq0}), and established the $L^2$ convergence in the infinite dimensional framework. Unfortunately, they deleted the distribution of the slow component in the fast equation. Another improvement is that Xu, Liu, Liu and Miao \cite{XLLM} inserted the distribution of the fast component into the fast part in the system (\ref{Eq0}), and also prove the $L^2$ convergence. 

In this paper, we study multiscale McKean-Vlasov stochastic systems where the whole systems depend on distributions of fast components, and establish an average principle in the $L^{2 p}(p \geqslant 1)$ sense. Concretely speaking, consider the following slow-fast system on $\mathbb{R}^n \times \mathbb{R}^m$ :
\be\left\{\begin{array}{l}
\dif X_{t}^{\e}=b_{1}(X_{t}^{\e},\sL^{\mP}_{X_{t}^{\e}},Z_{t}^{\e,z_0,\sL_\xi},\sL_{Z_{t}^{\e,\xi}})\dif t+\s_{1}(X_{t}^{\e},\sL^{\mP}_{X_{t}^{\e}})\dif B_{t},\\
X_{0}^{\e}=\varrho,\quad  0\leq t\leq T,\\
\dif Z_{t}^{\e,\xi}=\frac{1}{\e}b_{2}(X_{t}^{\e},\sL^{\mP}_{X_{t}^{\e}},Z_{t}^{\e,\xi},\sL^{\mP}_{Z_{t}^{\e,\xi}})\dif t+\frac{1}{\sqrt{\e}}\s_{2}(X_{t}^{\e},\sL^{\mP}_{X_{t}^{\e}},Z_{t}^{\e,\xi},\sL^{\mP}_{Z_{t}^{\e,\xi}})\dif W_{t},\\
Z_{0}^{\e,\xi}=\xi,\quad  0\leq t\leq T,\\
\dif Z_{t}^{\e,z_0,\sL^{\mP}_{\xi}}=\frac{1}{\e}b_{2}(X_{t}^{\e},\sL^{\mP}_{X_{t}^{\e}},Z_{t}^{\e,z_0,\sL^{\mP}_{\xi}},\sL^{\mP}_{Z_{t}^{\e,\xi}})\dif t+\frac{1}{\sqrt{\e}}\s_{2}(X_{t}^{\e},\sL^{\mP}_{X_{t}^{\e}},Z_{t}^{\e,z_0,\sL^{\mP}_{\xi}},\sL^{\mP}_{Z_{t}^{\e,\xi}})\dif W_{t},\\
Z_{0}^{\e,z_0,\sL^{\mP}_{\xi}}=z_0,\quad  0\leq t\leq T,
\end{array}
\right.
\label{Eq1}
\ee
where these mappings $b_1: \mathbb{R}^n \times \mathcal{P}_2\left(\mathbb{R}^n\right) \times \mathbb{R}^m \times \mathcal{P}_2\left(\mathbb{R}^m\right) \rightarrow \mathbb{R}^n, \sigma_1: \mathbb{R}^n \times \mathcal{P}_2\left(\mathbb{R}^n\right) \rightarrow \mathbb{R}^{n \times n}$, $b_2: \mathbb{R}^n \times \mathcal{P}_2\left(\mathbb{R}^n\right) \times \mathbb{R}^m \times \mathcal{P}_2\left(\mathbb{R}^m\right) \rightarrow \mathbb{R}^m, \sigma_2: \mathbb{R}^n \times \mathcal{P}_2\left(\mathbb{R}^n\right) \times \mathbb{R}^m \times \mathcal{P}_2\left(\mathbb{R}^m\right) \rightarrow \mathbb{R}^{m \times m}$ are all Borel measurable, and $\varrho, \xi$ are two random variables. To conclude the average principle for the system (\ref{Eq1}), since two frozen equations are McKean-Vlasov SDEs, their solutions are nonlinear Markov processes, and whether the classical Khasminskii time discretization or the Poisson equation to prove the strong convergence do not work. Therefore, we first linearize the two nonlinear Markov processes (cf. \cite{rrw}) and show the strong convergence by the modified Khasminskii time discretization method. Moreover, on account of these distributions, we apply a lot of tricks to obtain some estimations.

Next, nonlinear filtering problems mean that people extract some useful information of unobservable phenomena from observable ones, and estimate and predict them. So, the nonlinear filtering theory plays an important role in many areas including stochastic control, financial modeling, speech and image processing, and Bayesian networks (\cite{Goggin1, Goggin2, Imkeller, Qiao1, Qiao2, Qiao3}). And the nonlinear filtering theory of multiscale SDEs is systematically introduced by Kushner in \cite{kus}. Later more and 
more results about the nonlinear filtering of multiscale SDEs appear (See \cite{Goggin1, Goggin2, Imkeller, Qiao1, Qiao2, Qiao3} and references therein). However there are few results about nonlinear filtering of multiscale McKean-Vlasov SDEs. Hence, we also study the nonlinear filtering problem of them. That is, we define an observation process $Y_{t}^{\e}$ as follows
\be
Y_{t}^{\e}=V_{t}+\int_{0}^{t}h(X_{s}^{\e},\sL^{\mP}_{X_{s}^{\e}})\dif s,
\label{Eq4}
\ee
where $V_{\cdot}$ is a $l$-dimensional Brownian motion independent of $B_{\cdot}, W_{\cdot}$ and $h: \mR^n\times\cP_2(\mR^n)\mapsto\mR^l$ is Borel measurable. Then based on the obtained average principle, we establish that the nonlinear filtering of the slow part and its distribution converges to that of the average system in the $L^{q}$ ($p\geq 8, 1\leq q\leq \frac{p}{8}$) sense.

The novelty of this paper lies in three folds. The first fold is that the system (\ref{Eq1}) is more general than that in some known results (cf. \cite{ghl, hll, hlls, SunX, XLLM}). Thus, our result can be applied to many models. The second fold is that we obtain the strong convergence order for the $L^2$ case, which is important for numerical simulation. The third fold is that we prove the $L^{q}$ convergence of nonlinear filtering for multiscale McKean-Vlasov SDEs, which can partly cover some results in \cite{Imkeller, Qiao1, Qiao3}.

Lastly, we describe our motivation of this paper. Note that in \cite{XLLM}, although the multiscale system is general, the average principle is not right. This is because four authors used the Markov property which does not hold for general McKean-Vlasov SDEs. Our first motivation is to correct this mistake. Besides, as far as we know, no average principle for McKean-Vlasov SDEs with two time scales has yet been presented in the $L^{2p}$ ($p\geq 1$) sense. However, people usually need to estimate the higher order moments which possess a good robustness and can be applied in statistics, game theory, finance and other fields. So, it is our second motivation to establish an average principle in the $L^{2p}$ ($p\geq 1$) sense.

The paper proceeds as follows. In Section \ref{pre}, we introduce some related notations. Then we state main results in Section \ref{mare}. The proofs of two main theorems are placed in Section \ref{proofirs} and \ref{prooseco}, respectively.

The following convention will be used throughout the paper: $C$ with or without indices will denote different positive constants whose values may change from one place to another.

\section{Notations and assumptions}\label{pre}

In this section, we will recall some notations and list all assumptions.

\subsection{Notations}\label{nn}
In this subsection, we introduce some notations used in the sequel.

Let $|\cdot|, \|\cdot\|$ be the norms of a vector and a matrix, respectively. Let $\langle\cdot,\cdot\rangle$ be the inner product of vectors on $\mR^n$. $A^{*}$ denotes the transpose of the matrix $A$.

Let $\cB_b(\mR^{n})$ be the set of all bounded Borel measurable functions on $\mR^n$. Let $C(\mR^n)$ be the set of all  functions which are continuous on $\mR^n$. $C^{2}(\mR^n)$ represents the collection of all functions in $C(\mR^n)$ with continuous derivatives of order up to 2. 

Let $\sB(\mR^n)$ be the Borel $\sigma$-field on $\mR^n$. Let $\cP(\mR^n)$ be the collection of all probability measures on $\sB(\mR^n)$ with the usual topology of weak convergence. Let $\cP_{2}(\mR^n)$ denote  the collection of probability measures on $\sB(\mR^n)$ satisfying:
$$
\|\mu\|^{2}:=\int_{\mR^n}|x|^{2}\mu(dx)<\infty.
$$
It is known that $\cP_2(\mR^n)$ is a Polish space endowed with the $L^2$-Wasserstein distance defined by
$$
\mathbb{W}_2(\mu,\nu):= \inf\limits_{\pi\in\Psi(\mu,\nu)}\left(\int_{\mathbb{R}^n\times\mathbb{R}^n}|x-y|^{2}\pi(\dif x,\dif y)\right)^{\frac{1}{2}}, \quad \mu , \nu\in \cP_2(\mR^n),
$$
where $\Psi(\mu,\nu)$ is the set of all couplings $\pi$ with marginal distributions $\mu$ and $\nu$. Moreover, if $\xi,\zeta$ are two random variables with distributions $\sL_\xi, \sL_\zeta$ under $\mP$, respectively,
$$
\mathbb{W}_2(\sL_\xi, \sL_\zeta)\leq (\mE|\xi-\zeta|^2)^{\frac{1}{2}},
$$
where $\mE$ stands for the expectation with respect to $\mP$.

\subsection{Assumptions}\label{ass}

In this subsection, we give out all the assumptions used in the sequel:

\begin{enumerate}[$(\mathbf{H}^1_{b_{1}, \s_{1}})$]
\item
There exists a constant $L_{b_{1}, \s_{1}}>0$ such that for $x_{i}\in\mR^n$, $\mu_{i}\in\cP_{2}(\mR^n)$, $z_{i}\in\mR^m$, $\nu_i\in\cP_{2}(\mR^m)$, $i=1, 2$,
\ce
&&|b_{1}(x_{1},\mu_{1},z_{1},\nu_1)-b_{1}(x_{2},\mu_{2},z_{2},\nu_2)|^{2}+\|\s_{1}(x_{1},\mu_{1})-\s_{1}(x_{2},\mu_{2})\|^{2}\\
&\leq& L_{b_{1},\s_{1}}\(|x_{1}-x_{2}|^{2}+\mW_2^{2}(\mu_{1},\mu_{2})+|z_{1}-z_{2}|^{2}+\mW_2^{2}(\nu_{1},\nu_{2})\).
\de
\end{enumerate}
\begin{enumerate}[$(\mathbf{H}^1_{b_{2}, \s_{2}})$]
\item
There exists a constant $L_{b_{2}, \s_{2}}>0$ such that for $x_{i}\in\mR^n$, $\mu_{i}\in\cP_{2}(\mR^n)$, $z_{i}\in\mR^m$, $\nu_{i}\in\cP_{2}(\mR^m)$, $i=1, 2$,
\ce
&&|b_{2}(x_{1},\mu_{1},z_{1},\nu_1)-b_{2}(x_{2},\mu_{2},z_{2},\nu_2)|^{2}+\|\s_{2}(x_{1},\mu_{1},z_{1},\nu_1)-\s_{2}(x_{2},\mu_{2},z_{2},\nu_2)\|^{2}\\
&\leq& L_{b_{2}, \s_{2}}\(|x_{1}-x_{2}|^{2}+\mW_2^{2}(\mu_{1},\mu_{2})+|z_{1}-z_{2}|^{2}+\mW_2^{2}(\nu_{1},\nu_{2})\).
\de
\end{enumerate}
\begin{enumerate}[$(\mathbf{H}^2_{b_{2}, \s_{2}})$]
\item
For some $p\geq 0$, there exist two constants $\b_1>0, \b_2>0$ satisfying $\b_{1}-\b_{2}>(4p+4)L_{b_2,\sigma_2}$ such that for $x\in\mR^n$, $\mu\in\cP_{2}(\mR^n)$, $z_{i}\in\mR^m$, $\nu_{i}\in\cP_{2}(\mR^m)$, $i=1, 2$,
\ce
&&2\<z_{1}-z_{2},b_{2}(x,\mu,z_{1},\nu_{1})-b_{2}(x,\mu,z_{2},\nu_{2})\>
+(2p+1)\|\s_{2}(x,\mu,z_{1},\nu_{1})-\s_{2}(x,\mu,z_{2},\nu_{2})\|^{2}\\
&\leq& -\b_1|z_{1}-z_{2}|^{2}+\b_2\mW_2^{2}(\nu_{1},\nu_{2}).
\de
\end{enumerate}
\begin{enumerate}[$(\mathbf{H}_{h})$]
\item $h$ is bounded, and there is a constant $L_{h}>0$ such that
$$
|h(x_{1},\mu_{1})-h(x_{2},\mu_{2})|^{2}\leq L_{h}(|x_{1}-x_{2}|^{2}+\mW_2^{2}(\mu_{1},\mu_{2})).
$$
\end{enumerate}

\br
(i) $(\mathbf{H}^1_{b_{1}, \s_{1}})$ yields that there exists a constant $\bar{L}_{b_{1}, \s_{1}}>0$ such that for $x\in\mR^n$, $\mu\in\cP_{2}(\mR^n)$, $z\in\mR^m$, $\nu\in\cP_{2}(\mR^m)$
\be
|b_{1}(x,\mu,z,\nu)|^{2}+\|\s_{1}(x,\mu)\|^{2}\leq \bar{L}_{b_{1}, \s_{1}}(1+|x|^{2}+\|\mu\|^{2}+|z|^{2}+\|\nu\|^{2}).
\label{b1line}
\ee

(ii) $(\mathbf{H}^1_{b_{2}, \s_{2}})$ implies that there exists a constant $\bar{L}_{b_{2}, \s_{2}}>0$ such that for $x\in\mR^n$, $\mu\in\cP_{2}(\mR^n)$, $z\in\mR^m$, $\nu\in\cP_{2}(\mR^m)$,
\be
|b_{2}(x,\mu,z,\nu)|^{2}+\|\s_{2}(x,\mu,z,\nu)\|^{2}
\leq \bar{L}_{b_{2}, \s_{2}}(1+|x|^{2}+\|\mu\|^{2}+|z|^{2}+\|\nu\|^{2}).
\label{b2nu}
\ee

(iii) By $(\mathbf{H}^1_{b_{2}, \s_{2}})$ and $(\mathbf{H}^2_{b_{2}, \s_{2}})$, it holds that for $x\in\mR^n$, $\mu\in\cP_{2}(\mR^n)$, $z\in\mR^m$, $\nu\in\cP_{2}(\mR^m)$
\be
2\<z,b_{2}(x,\mu,z,\nu)\>+(2p+1)\|\s_{2}(x,\mu,z,\nu)\|^{2}\leq -\a_{1}|z|^{2}+\a_{2}\|\nu\|^{2}+C(1+|x|^{2}+\|\mu\|^{2}),
\label{bemu}
\ee
where $\a_{1}:=\b_1-(2p+2)L_{b_{2}, \s_{2}}$, $\a_{2}:=\b_{2}+(2p+1)L_{b_{2}, \s_{2}}$, $\a_1-\a_2-L_{b_{2}, \s_{2}}>0$ and $C>0$ is a constant.
\er

\section{Main results}\label{mare}

In this section, we provide main results of the paper. 

\subsection{The average principle for multiscale McKean-Vlasov SDEs}

In this subsection, we state the average principle result for multiscale McKean-Vlasov SDEs.

Let us recall the system (\ref{Eq1}), i.e.
\ce\left\{\begin{array}{l}
\dif X_{t}^{\e}=b_{1}(X_{t}^{\e},\sL^{\mP}_{X_{t}^{\e}},Z_{t}^{\e,z_0,\sL_\xi},\sL_{Z_{t}^{\e,\xi}})\dif t+\s_{1}(X_{t}^{\e},\sL^{\mP}_{X_{t}^{\e}})\dif B_{t},\\
X_{0}^{\e}=\varrho,\quad  0\leq t\leq T,\\
\dif Z_{t}^{\e,\xi}=\frac{1}{\e}b_{2}(X_{t}^{\e},\sL^{\mP}_{X_{t}^{\e}},Z_{t}^{\e,\xi},\sL^{\mP}_{Z_{t}^{\e,\xi}})\dif t+\frac{1}{\sqrt{\e}}\s_{2}(X_{t}^{\e},\sL^{\mP}_{X_{t}^{\e}},Z_{t}^{\e,\xi},\sL^{\mP}_{Z_{t}^{\e,\xi}})\dif W_{t},\\
Z_{0}^{\e,\xi}=\xi,\quad  0\leq t\leq T,\\
\dif Z_{t}^{\e,z_0,\sL^{\mP}_{\xi}}=\frac{1}{\e}b_{2}(X_{t}^{\e},\sL^{\mP}_{X_{t}^{\e}},Z_{t}^{\e,z_0,\sL^{\mP}_{\xi}},\sL^{\mP}_{Z_{t}^{\e,\xi}})\dif t+\frac{1}{\sqrt{\e}}\s_{2}(X_{t}^{\e},\sL^{\mP}_{X_{t}^{\e}},Z_{t}^{\e,z_0,\sL^{\mP}_{\xi}},\sL^{\mP}_{Z_{t}^{\e,\xi}})\dif W_{t},\\
Z_{0}^{\e,z_0,\sL^{\mP}_{\xi}}=z_0,\quad  0\leq t\leq T,
\end{array}
\right.
\de
where $\mathbb{E}|\varrho|^{2 p+2}<\infty, \mathbb{E}|\xi|^{2 p+2}<\infty$ ($p$ is the same to that in $(\mathbf{H}_{b_2,\sigma_2}^2)$).
Under $(\mathbf{H}^1_{b_{1}, \s_{1}})$ $(\mathbf{H}^1_{b_{2}, \s_{2}})$, by \cite[Theorem 2.1]{WangFY}, the system (\ref{Eq1}) has a unique strong solution $(X_{\cdot}^{\e}, Z_{\cdot}^{\e,\xi},Z_{\cdot}^{\e,z_0,\sL^{\mP}_{\xi}})$.

Next, we take any $x\in \mR^{n}$ and $\mu\in\cP_{2}(\mR^n)$, and fix them. Consider the following SDE:
\be\left\{\begin{array}{l}
\dif Z_{t}^{x,\mu,\xi}=b_{2}(x,\mu,Z_{t}^{x,\mu,\xi},\sL^{\mP}_{Z_{t}^{x,\mu,\xi}})\dif t+\s_{2}(x,\mu,Z_{t}^{x,\mu,\xi},\sL^{\mP}_{Z_{t}^{x,\mu,\xi}})\dif W_{t},\\
Z_{0}^{x,\mu,\xi}=\xi, \quad 0 \leq t \leq T.
\end{array}
\right.
\label{Eq2}
\ee
Under $(\mathbf{H}^1_{b_{2}, \s_{2}})$, by \cite[Theorem 2.1]{WangFY} we know that the above equation has a unique strong solution $Z_{\cdot}^{x,\mu,\xi}$. Moreover, under $(\mathbf{H}^2_{b_{2}, \s_{2}})$, by \cite[Theorem 3.1]{WangFY}, one could obtain that there exists a unique invariant probability measure $\eta^{x,\mu}$ for Eq.(\ref{Eq2}). So, we construct an average equation on $(\Omega,\sF,\{\sF_{t}\}_{t\in[0,T]},\mP)$ as follows:
\be\left\{\begin{array}{l}
\dif \bar{X}_{t}=\bar{b}_{1}(\bar{X}_{t},\sL^{\mP}_{\bar{X}_{t}})\dif t+\s_1(\bar{X}_{t},\sL^{\mP}_{\bar{X}_{t}})\dif B_{t},\\
\bar{X}_{0}=\varrho,
\end{array}
\right.
\label{Eq3}
\ee
where $\bar{b}_{1}(x,\mu)=\int_{\mR^{m}\times\cP_2(\mR^m)}b_{1}(x,\mu,z,\nu)\eta^{x,\mu}\times\d_{\eta^{x,\mu}}(\dif z,\dif \nu)$. 

Now, it is the position to state the first main result.

\bt \label{xbarxp}
Under these assumptions $(\mathbf{H}^1_{b_{1}, \s_{1}})$ $(\mathbf{H}^{1}_{b_{2}, \s_{2}})$-$(\mathbf{H}^{2}_{b_{2}, \s_{2}})$, for $p\geq 1$, it holds that
\be
\lim\limits_{\e\rightarrow 0}\mE\(\sup_{0\leq t\leq T}|X_{t}^{\e}-\bar{X}_{t}|^{2p}\)=0,
\label{mesu}
\ee
where $\bar{X}$ is a solution of Eq.(\ref{Eq3}). In particular, we have that for any $0<\g<1$
\be
\mE\(\sup_{0\leq t\leq T}|X_{t}^{\e}-\bar{X}_{t}|^{2}\)\leq C(\e^{1-\g}+\e^{2\g}+\e^\g).
\label{2or}
\ee
\et

The proof of Theorem \ref{xbarxp} is placed in Section \ref{proofirs}.

\br
For (\ref{2or}), if we take $\g=1/2$, it follows that
\ce
\mE\(\sup_{0\leq t\leq T}|X_{t}^{\e}-\bar{X}_{t}|^{2}\)\leq C\e^{1/2}.
\de
That is, we obtain the convergence order $1/4$.
\er

\subsection{The efficient filtering for multiscale McKean-Vlasov SDEs}

In this subsection, we state the efficient filtering result for multiscale McKean-Vlasov SDEs.

Set
\ce
(\Lambda_{t}^{\e})^{-1}:=\exp\left\{-\int_{0}^{t}h^{i}(X_{s}^{\e},\sL^{\mP}_{X_{s}^{\e}})\dif V_{s}^{i} -\frac{1}{2}\int_{0}^{t}|h(X_{s}^{\e},\sL^{\mP}_{X_{s}^{\e}})|^{2}\dif s\right\}.
\de
Here and hereafter, we use the convention that repeated indices imply summation. Under $(\mathbf{H}_{h})$, we get that
$$
\mE\left(\exp\left\{\frac{1}{2}\int_{0}^{T}|h(X_{s}^{\e},\sL^{\mP}_{X_{s}^{\e}})|^{2}\dif s\right\}\right)<\infty,
$$
and furthermore $(\Lambda_{t}^{\e})^{-1}$ is an exponential martingale under the measure $\mP$. Define a probability measure $\mP^{\e}$ via
\ce
\frac{\dif\mP^{\e}}{\dif \mP}=(\Lambda_{T}^{\e})^{-1}.
\de
Then by the Girsanov theorem, it holds that $Y_{\cdot}^{\e}$ is a Brownian motion under the probability measure $\mP^{\e}$.

Define the nonlinear filtering for the state $X_{t}^{\e}$ and the measure $\sL^{\mP}_{X_{t}^{\e}}$:
\ce
\rho_{t}^{\e}(F)
&:=&\mE^{\mP^{\e}}[F(X_{t}^{\e},\sL^{\mP}_{X_{t}^{\e}})\Lambda_{t}^{\e}|\mathscr{F}_{t}^{Y^{\e}}],\\
\pi_{t}^{\e}(F)
&:=&\mE[F(X_{t}^{\e},\sL^{\mP}_{X_{t}^{\e}})|\mathscr{F}_{t}^{Y^{\e}}],\quad F\in\cB_b(\mR^n\times\cP_2(\mR^n)),
\de
where $\mathscr{F}_{t}^{Y^{\e}}=\sigma\{Y_{s}^{\e},0\leq s \leq t\}\vee \cN$, and $\cN$ denotes  the collection of all zero sets under $\mP$. Here $\rho_{t}^{\e}(F)$, $\pi_{t}^{\e}(F)$ are called the unnormalized filtering and the normalized filtering of $(X_{t}^{\e},  \sL^{\mP}_{X_{t}^{\e}})$  with respect to $\mathscr{F}_{t}^{Y^{\e}}$, respectively. By the Kallianpur-Striebel formula, we get the following relationship between $\rho_{t}^{\e}(F)$ and $\pi_{t}^{\e}(F)$:
$$
\pi_{t}^{\e}(F)=\frac{\rho_{t}^{\e}(F)}{\rho_{t}^{\e}(1)}.
$$

Next, set
\ce
&&\Lambda_{t}^{0}:=\exp\left\{\int_{0}^{t}h^{i}(\bar{X}_{s},\sL^{\mP}_{\bar{X}_{s}})\dif Y_{s}^{\e, i} -\frac{1}{2}\int_{0}^{t}|h(\bar{X}_{s},\sL^{\mP}_{\bar{X}_{s}})|^{2}\dif s\right\},\\
&&\rho_{t}^{0}(F):=\mE^{\mP^{\e}}[F(\bar{X}_{t},\sL^{\mP}_{\bar{X}_{t}})\Lambda_{t}^{0}|\mathscr{F}_{t}^{Y^{\varepsilon}}],\\
&&\pi_{t}^{0}(F):=\frac{\rho_{t}^{0}(F)}{\rho_{t}^{0}(1)},
\de
and we study the relationship between $\pi_{t}^{\e}$ and $\pi_{t}^{0}$. The second main result of the paper is the following theorem.

\bt\label{effifilt}
Under these assumptions $(\mathbf{H}^1_{b_{1}, \s_{1}})$, $(\mathbf{H}^{1}_{b_{2}, \s_{2}})$-$(\mathbf{H}^{2}_{b_{2}, \s_{2}})$ and $(\mathbf{H}_{h})$, for $p\geq 8$, the nonlinear filtering of the original system converges to that of the average system in the $L^{q}$ ($1\leq q\leq \frac{p}{8}$) sense under $\mP$, that is,
\be
\lim_{\e\rightarrow 0}\mE|\pi_{t}^{\e}(F)-\pi_{t}^{0}(F)|^q=0, \quad F\in C_{b,lip}(\mR^n\times\cP_2(\mR^n)),
\label{8}
\ee
where $C_{b,lip}(\mR^n\times\cP_2(\mR^n))$ stands for the collection of all bounded and Lipschitz continuous functions on $\mR^n\times\cP_2(\mR^n)$.
\et

The proof of Theorem \ref{effifilt} is placed in Section \ref{prooseco}.

\section{Proof of Theorem \ref{xbarxp}}\label{proofirs}

In the section, we prove Theorem \ref{xbarxp}. The proof consists of three parts. In the first part (Subsection \ref{xtzthatxhatz}), we segment the time interval $[0, T]$ by the size $\d$, where $\d$ is a fixed positive number depending on $\e$, and introduce three auxiliary processes: 
\be
&&\hat{Z}_{t}^{\e,\xi}
=\xi+\frac{1}{\e}\int_{0}^tb_{2}(X_{s(\d)}^{\e},\sL^{\mP}_{X_{s(\d)}^{\e}},\hat{Z}_{s}^{\e,\xi},\sL^{\mP}_{\hat{Z}_{s}^{\e,\xi}})\dif s
+\frac{1}{\sqrt{\e}}\int_{0}^t\s_{2}(X_{s(\d)}^{\e},\sL^{\mP}_{X_{s(\d)}^{\e}},\hat{Z}_{s}^{\e,\xi},\sL^{\mP}_{\hat{Z}_{s}^{\e,\xi}})\dif W_{s},\no\\
\label{hatz}\\
&&\hat{Z}_{t}^{\e,z_0,\sL^{\mP}_{\xi}}
=z_0+\frac{1}{\e}\int_{0}^tb_{2}(X_{s(\d)}^{\e},\sL^{\mP}_{X_{s(\d)}^{\e}},\hat{Z}_{s}^{\e,z_0,\sL^{\mP}_{\xi}},\sL^{\mP}_{\hat{Z}_{s}^{\e,\xi}})\dif s\no\\
&&\qquad\qquad\qquad+\frac{1}{\sqrt{\e}}\int_{0}^t\s_{2}(X_{s(\d)}^{\e},\sL^{\mP}_{X_{s(\d)}^{\e}},\hat{Z}_{s}^{\e,z_0,\sL^{\mP}_{\xi}},\sL^{\mP}_{\hat{Z}_{s}^{\e,\xi}})\dif W_{s},\label{hatz2}\\
&&\hat{X}_t^\e=\varrho+\int_0^tb_1(X_{s(\d)}^{\e},\sL^{\mP}_{X_{s(\d)}^{\e}},\hat{Z}_{s}^{\e,z_0,\sL^{\mP}_{\xi}},\sL^{\mP}_{\hat{Z}_{s}^{\e,\xi}})\dif s+\int_0^t\s_{1}(X_{s}^{\e},\sL^{\mP}_{X_{s}^{\e}})\dif B_{s}, \label{hatx}
\ee
where $s(\d)=[\frac{s}{\d}]\d$, and $[\frac{s}{\d}]$ denotes the integer part of $\frac{s}{\d}$. Then we estimate $X_{\cdot}^{\e}, Z_{\cdot}^{\e,\xi},Z_{\cdot}^{\e,z_0,\sL^{\mP}_{\xi}}$, $\hat{X}_{\cdot}^{\e}, \hat{Z}_{\cdot}^{\e,\xi},\hat{Z}_{\cdot}^{\e,z_0,\sL^{\mP}_{\xi}}$. In the second part (Subsection \ref{froequ}) and the third part (Subsection \ref{aveequ}), we present some estimates for the frozen equation (\ref{Eq2}) and the average equation (\ref{Eq3}), respectively.

\subsection{Some estimates for $X_{\cdot}^{\e}, Z_{\cdot}^{\e,\xi},Z_{\cdot}^{\e,z_0,\sL^{\mP}_{\xi}},$ $\hat{X}_{\cdot}^{\e}, \hat{Z}_{\cdot}^{\e,\xi},\hat{Z}_{\cdot}^{\e,z_0,\sL^{\mP}_{\xi}}$}\label{xtzthatxhatz}

\bl \label{xtztc}
 Under assumptions $(\mathbf{H}^{1}_{b_{1}, \s_{1}})$ $(\mathbf{H}^{1}_{b_{2}, \s_{2}})$ $(\mathbf{H}^{2}_{b_{2}, \s_{2}})$, there exists a constant $C>0$ such that 
\ce
&&\sup\limits_{\e}\mE\left(\sup\limits_{t\in[0,T]}|X_{t}^{\e}|^{2p+2}\right)\leq C(1+\mE|\varrho|^{2p+2}+|z_0|^{2p+2}+\mE|\xi|^{2p+2}),\\
&&\sup\limits_{t\in[0,T]}\mE|Z_{t}^{\e,\xi}|^{2p+2}\leq C(1+\mE|\varrho|^{2p+2}+|z_0|^{2p+2}+\mE|\xi|^{2p+2}),\\
&&\sup\limits_{t\in[0,T]}\mE|Z_{t}^{\e,z_0,\sL^{\mP}_{\xi}}|^{2p+2}\leq C(1+\mE|\varrho|^{2p+2}+|z_0|^{2p+2}+\mE|\xi|^{2p+2}).
\de
\el
\begin{proof}
For $X_{t}^{\e}$, based on the BDG inequality and $(\mathbf{H}^{1}_{b_{1}, \s_{1}})$, we can get
\be
\mE\left(\sup\limits_{s\in[0,t]}|X_{s}^{\e}|^{2p+2}\right)
&\leq& 3^{2p+1}\mE|\varrho|^{2p+2}+3^{2p+1}\mE\left(\sup\limits_{s\in[0,t]}\Big|\int_{0}^{s}b_{1}(X_{r}^{\e},\sL^{\mP}_{X_{r}^{\e}},Z_{r}^{\e,z_0,\sL^{\mP}_{\xi}},\sL^{\mP}_{Z_{r}^{\e,\xi}})\dif r\Big|^{2p+2}\right)\no\\
&&+3^{2p+1}\mE\left(\sup\limits_{s\in[0,t]}\Big|\int_{0}^{s}\s_{1}(X_{r}^{\e},\sL^{\mP}_{X_{r}^{\e}})\dif B_{r}\Big|^{2p+2}\right)\no\\
&\leq& 3^{2p+1}\mE|\varrho|^{2p+2}+(3t)^{2p+1}\mE\int_{0}^{t}\Big|b_{1}(X_{r}^{\e},\sL^{\mP}_{X_{r}^{\e}},Z_{r}^{\e,z_0,\sL^{\mP}_{\xi}},\sL^{\mP}_{Z_{r}^{\e,\xi}})\Big|^{2p+2}\dif r\no\\
&&+3^{2p+1}t^{p}C\mE\int_{0}^{t}\|\s_{1}(X_{r}^{\e},\sL^{\mP}_{X_{r}^{\e}})\|^{2p+2}\dif r\no\\
&\leq&3^{2p+1}\mE|\varrho|^{2p+2}+C\mE\int_{0}^{t}(1+|X_{r}^{\e}|+\|\sL^{\mP}_{X_{r}^{\e}}\|+|Z_{r}^{\e,z_0,\sL^{\mP}_{\xi}}|+\|\sL^{\mP}_{Z_{r}^{\e,\xi}}\|)^{2p+2}\dif r\no\\
&\leq&C(\mE|\varrho|^{2p+2}+1)+C\int_{0}^{t}\mE|X_{r}^{\e}|^{2p+2}\dif r+C\int_{0}^{t}\mE|Z_{r}^{\e,z_0,\sL^{\mP}_{\xi}}|^{2p+2}\dif r\no\\
&&+C\int_{0}^{t}\mE|Z_{r}^{\e,\xi}|^{2p+2}\dif r,
\label{exqc}
\ee
where $\|\sL^{\mP}_{X_{r}^{\e}}\|^2=\mE|X_{r}^{\e}|^2$.

For $Z_{t}^{\e,\xi}$, applying the It\^{o} formula to $|Z_{t}^{\e,\xi}|^{2p+2}$ and taking the expectation, one could obtain that
\ce
\mE|Z_{t}^{\e,\xi}|^{2p+2}
&=& \mE|\xi|^{2p+2}+\frac{2p+2}{\e}\mE\int_{0}^{t}|Z_{s}^{\e,\xi}|^{2p}\<Z_{s}^{\e,\xi}, b_{2}(X_{s}^{\e},\sL^{\mP}_{X_{s}^{\e}},Z_{s}^{\e,\xi},\sL^{\mP}_{Z_{s}^{\e,\xi}})\>\dif s\\
&&+\frac{2p(p+1)}{\e}\mE\int_{0}^{t}|Z_{s}^{\e,\xi}|^{2p-2}\|\s_{2}(X_{s}^{\e},\sL^{\mP}_{X_{s}^{\e}},Z_{s}^{\e,\xi},\sL^{\mP}_{Z_{s}^{\e,\xi}})Z_{s}^{\e,\xi}\|^2\dif s\\
&&+\frac{p+1}{\e}\mE\int_{0}^{t}|Z_{s}^{\e,\xi}|^{2p}\|\s_{2}(X_{s}^{\e},\sL^{\mP}_{X_{s}^{\e}},Z_{s}^{\e,\xi},\sL^{\mP}_{Z_{s}^{\e,\xi}})\|^2\dif s,
\de
and 
\ce
\frac{\dif}{\dif t}\mE|Z_{t}^{\e,\xi}|^{2p+2}&=&\frac{2p+2}{\e}\mE|Z_{t}^{\e,\xi}|^{2p}\<Z_{t}^{\e,\xi}, b_{2}(X_{t}^{\e},\sL^{\mP}_{X_{t}^{\e}},Z_{t}^{\e,\xi},\sL^{\mP}_{Z_{t}^{\e,\xi}})\>\\
&&+\frac{2p(p+1)}{\e}\mE|Z_{t}^{\e,\xi}|^{2p-2}\|\s_{2}(X_{t}^{\e},\sL^{\mP}_{X_{t}^{\e}},Z_{t}^{\e,\xi},\sL^{\mP}_{Z_{t}^{\e,\xi}})Z_{t}^{\e,\xi}\|^2\\
&&+\frac{p+1}{\e}\mE|Z_{t}^{\e,\xi}|^{2p}\|\s_{2}(X_{t}^{\e},\sL^{\mP}_{X_{t}^{\e}},Z_{t}^{\e,\xi},\sL^{\mP}_{Z_{t}^{\e,\xi}})\|^2\\
&\leq&\frac{2p+2}{\e}\mE|Z_{t}^{\e,\xi}|^{2p}\<Z_{t}^{\e,\xi}, b_{2}(X_{t}^{\e},\sL^{\mP}_{X_{t}^{\e}},Z_{t}^{\e,\xi},\sL^{\mP}_{Z_{t}^{\e,\xi}})\>\\
&&+\frac{(2p+1)(p+1)}{\e}\mE|Z_{t}^{\e,\xi}|^{2p}\|\s_{2}(X_{t}^{\e},\sL^{\mP}_{X_{t}^{\e}},Z_{t}^{\e,\xi},\sL^{\mP}_{Z_{t}^{\e,\xi}})\|^2\\
&\leq&\frac{p+1}{\e}\mE|Z_{t}^{\e,\xi}|^{2p}\(-\a_{1}|Z_{t}^{\e,\xi}|^{2}+\a_{2}\|\sL^{\mP}_{Z_{t}^{\e,\xi}}\|^{2}+C(1+|X_{t}^{\e}|^{2}+\|\sL^{\mP}_{X_{t}^{\e}}\|^{2})\)\\
&\leq&\frac{p+1}{\e}\bigg[-\a_{1}\mE|Z_{t}^{\e,\xi}|^{2p+2}+\a_{2}\mE|Z_{t}^{\e,\xi}|^{2p+2}+L_{b_2,\s_2}\mE|Z_{t}^{\e,\xi}|^{2p+2}\\
&&\quad\qquad+C(1+\mE|X_{t}^{\e}|^{2p+2}+\|\sL^{\mP}_{X_{t}^{\e}}\|^{2p+2})\bigg]\\
&\leq&\frac{-(\a_1-\a_2-L_{b_2,\s_2})(p+1)}{\e}\mE|Z_{t}^{\e,\xi}|^{2p+2}+\frac{C}{\e}(\mE|X_{t}^{\e}|^{2p+2}+1),
\de
where the above inequality is based on (\ref{bemu}). By the comparison theorem, we have that
\be
\mE|Z_{t}^{\e,\xi}|^{2p+2}&\leq&\mE|\xi|^{2p+2}e^{-\frac{(\a_1-\a_2-L_{b_2,\s_2})(p+1)}{\e}t}+\frac{C}{\e}\int_{0}^{t}e^{-\frac{(\a_1-\a_2-L_{b_2,\s_2})(p+1)}{\e}(t-s)}(\mE|X_{s}^{\e}|^{2p+2}+1)\dif s\no\\
&\leq&\mE|\xi|^{2p+2}+C\left(\mE\left(\sup\limits_{s\in[0,t]}|X_{s}^{\e}|^{2p+2}\right)+1\right).
\label{zees}
\ee

Next, for $Z_{\cdot}^{\e,z_0,\sL^{\mP}_{\xi}}$, by the similar deduction to that for $Z_{t}^{\e,\xi}$, it holds that
\be
\mE|Z_t^{\e,z_0,\sL^{\mP}_{\xi}}|^{2p+2}\leq |z_0|^{2p+2}+\mE|\xi|^{2p+2}+C\left(\mE\left(\sup\limits_{s\in[0,t]}|X_{s}^{\e}|^{2p+2}\right)+1\right).
\label{zezes}
\ee

Inserting (\ref{zees}) (\ref{zezes}) in (\ref{exqc}), by the Gronwall inequality one can get that
$$
\mE\left(\sup\limits_{t\in[0,T]}|X_{t}^{\e}|^{2p+2}\right)\leq C(1+\mE|\varrho|^{2p+2}+|z_0|^{2p+2}+\mE|\xi|^{2p+2}),
$$
which together with (\ref{zees}) (\ref{zezes}) implies that
\ce
&&\sup\limits_{t\in[0,T]}\mE|Z_{t}^{\e,\xi}|^{2p+2}\leq C(1+\mE|\varrho|^{2p+2}+|z_0|^{2p+2}+\mE|\xi|^{2p+2}),\\
&&\sup\limits_{t\in[0,T]}\mE|Z_{t}^{\e,z_0,\sL^{\mP}_{\xi}}|^{2p+2}\leq C(1+\mE|\varrho|^{2p+2}+|z_0|^{2p+2}+\mE|\xi|^{2p+2}).
\de
The proof is complete.
\end{proof}

Next, we estimate $\mE|X_{t}^{\e}-X_{k\d}^{\e}|^{2}$ for any $t\in[k\d, (k+1)\d)$ and $k=0,1,2,\cdots,[\frac{T}{\d}]-1$. Note that
\ce
X_{t}^{\e}-X_{k\d}^{\e}=\int_{k\d}^{t}b_{1}(X_{s}^{\e},\sL^{\mP}_{X_{s}^{\e}},Z_{s}^{\e,z_0,\sL_\xi},\sL_{Z_{s}^{\e,\xi}})\dif s
+\int_{k\d}^{t}\s_{1}(X_{s}^{\e},\sL^{\mP}_{X_{s}^{\e}})\dif B_{s}.
\de
By $(\ref{b1line})$ and the BDG inequality, it holds that
\be
&&\mE|X_{t}^{\e}-X_{k\d}^{\e}|^{2}\no\\
&\leq& 2\bigg(\mE\left|\int_{k\d}^{t}b_{1}(X_{s}^{\e},\sL^{\mP}_{X_{s}^{\e}},Z_{s}^{\e,z_0,\sL^{\mP}_\xi},\sL^{\mP}_{Z_{s}^{\e,\xi}})\dif s\right|^{2}
+\mE\left|\int_{k\d}^{t}\s_{1}(X_{s}^{\e},\sL^{\mP}_{X_{s}^{\e}})\dif B_{s}\right|^{2}\bigg)\no\\
&\leq& 2\bigg(\d\int_{k\d}^{t}\mE\left|b_{1}(X_{s}^{\e},\sL^{\mP}_{X_{s}^{\e}},Z_{s}^{\e,z_0,\sL^{\mP}_\xi},\sL^{\mP}_{Z_{s}^{\e,\xi}})\right|^{2}\dif s
+\int_{k\d}^{t}\mE\left\|\s_{1}(X_{s}^{\e},\sL^{\mP}_{X_{s}^{\e}})\right\|^{2}\dif s\bigg)\no\\
&\leq& C({\d}^{2}+{\d}),
\label{barx}
\ee
where the last inequality is based on Lemma \ref{xtztc}.

 Moreover, by the same deduction to that for $Z_{\cdot}^{\e,\xi}, Z_{\cdot}^{\e,z_0,\sL^{\mP}_\xi}$ in Lemma \ref{xtztc}, we obtain the following estimate.
\bl\label{hatze}
Under these assumptions $(\mathbf{H}^1_{b_{2}, \s_{2}})$-$(\mathbf{H}^2_{b_{2}, \s_{2}})$, it holds that 
\ce
&&\sup\limits_{t\in[0,T]}\mE|\hat{Z}_{t}^{\e,\xi}|^{2}\leq C(1+\mE|\varrho|^{2}+|z_0|^{2}+\mE|\xi|^{2}),\\
&&\sup\limits_{t\in[0,T]}\mE|\hat{Z}_{t}^{\e,z_0,\sL^{\mP}_{\xi}}|^{2}\leq C(1+\mE|\varrho|^{2}+|z_0|^{2}+\mE|\xi|^{2}).
\de
\el

\bl\label{zehatze}
 Suppose $(\mathbf{H}^{1}_{b_{1},\s_{1}})$, $(\mathbf{H}^{1}_{b_{2},\s_{2}})$, $(\mathbf{H}^{2}_{b_{2},\s_{2}})$ hold. Then for $p\geq 1$, there exists a constant $C>0$ such that 
 \ce
&&\sup\limits_{t\in[0,T]}\mE|Z_{t}^{\e,\xi}-\hat{Z}_{t}^{\e,\xi}|^{2}\leq \frac{C}{\b_1-\b_2-L_{b_{2},\s_2}}
(\d^2+\d),\\
&&\sup\limits_{t\in[0,T]}\mE|Z_{t}^{\e,z_0,\sL^{\mP}_\xi}-\hat{Z}_{t}^{\e,z_0,\sL^{\mP}_\xi}|^{2}\leq \frac{C}{\b_1-L_{b_{2},\s_2}}({\d}^{2}+{\d}).
\de
\el
\begin{proof}
First of all, by (\ref{Eq1}) and (\ref{hatz}), we have that 
\ce
Z_{t}^{\e,\xi}-\hat{Z}_{t}^{\e,\xi}
&=&\frac{1}{\e}\int_{0}^{t}\(b_{2}(X_{s}^{\e},\sL^{\mP}_{X_{s}^{\e}},Z_{s}^{\e,\xi},\sL^{\mP}_{Z_{s}^{\e,\xi}})
-b_{2}(X_{s(\d)}^{\e},\sL^{\mP}_{X_{s(\d)}^{\e}},\hat{Z}_{s}^{\e,\xi},\sL^{\mP}_{\hat{Z}_{s}^{\e,\xi}})\)\dif s\\
&&+\frac{1}{\sqrt{\e}}\int_{0}^{t}\(\s_{2}(X_{s}^{\e},\sL^{\mP}_{X_{s}^{\e}},Z_{s}^{\e,\xi},\sL^{\mP}_{Z_{s}^{\e,\xi}})
-\s_{2}(X_{s(\d)}^{\e},\sL^{\mP}_{X_{s(\d)}^{\e}},\hat{Z}_{s}^{\e,\xi},\sL^{\mP}_{\hat{Z}_{s}^{\e,\xi}})\)\dif W_{s}.
\de
Applying the It\^{o} formula to $|Z_{t}^{\e,\xi}-\hat{Z}_{t}^{\e,\xi}|^{2}$ and taking the expectation, by $(\mathbf{H}^{1}_{b_{2},\s_{2}})$ and $(\mathbf{H}^{2}_{b_{2},\s_{2}})$ one could obtain that
\ce
&&\mE|Z_{t}^{\e,\xi}-\hat{Z}_{t}^{\e,\xi}|^{2}\\
&=&\frac{1}{\e}
\mE\int_{0}^{t}2\<Z_{s}^{\e,\xi}-\hat{Z}_{s}^{\e,\xi}, b_{2}(X_{s}^{\e},\sL^{\mP}_{X_{s}^{\e}},Z_{s}^{\e,\xi},\sL^{\mP}_{Z_{s}^{\e,\xi}})
-b_{2}(X_{s(\d)}^{\e},\sL^{\mP}_{X_{s(\d)}^{\e}},\hat{Z}_{s}^{\e,\xi},\sL^{\mP}_{\hat{Z}_{s}^{\e,\xi}})\>\dif s\\
&&+\frac{1}{\e}
\mE\int_{0}^{t}\|\s_{2}(X_{s}^{\e},\sL^{\mP}_{X_{s}^{\e}},Z_{s}^{\e,\xi},\sL^{\mP}_{Z_{s}^{\e,\xi}})
-\s_{2}(X_{s(\d)}^{\e},\sL^{\mP}_{X_{s(\d)}^{\e}},\hat{Z}_{s}^{\e,\xi},\sL^{\mP}_{\hat{Z}_{s}^{\e,\xi}})\|^{2}\dif s
\de
and
\ce
&&\frac{\dif \mE|Z_{t}^{\e,\xi}-\hat{Z}_{t}^{\e,\xi}|^{2}}{\dif t}\\
 &\leq&\frac{1}{\e}
 \mE\bigg[2\<Z_{t}^{\e,\xi}-\hat{Z}_{t}^{\e,\xi}, b_{2}(X_{t}^{\e},\sL^{\mP}_{X_{t}^{\e}},Z_{t}^{\e,\xi},\sL^{\mP}_{Z_{t}^{\e,\xi}})
 -b_{2}(X_{t}^{\e},\sL^{\mP}_{X_{t}^{\e}},\hat{Z}_{t}^{\e,\xi},\sL^{\mP}_{\hat{Z}_{t}^{\e,\xi}})\>\\
 &&+(2p+1)\|\s_{2}(X_{t}^{\e},\sL^{\mP}_{X_{t}^{\e}},Z_{t}^{\e,\xi},\sL^{\mP}_{Z_{t}^{\e,\xi}})
 -\s_{2}(X_{t}^{\e},\sL^{\mP}_{X_{t}^{\e}},\hat{Z}_{t}^{\e,\xi},\sL^{\mP}_{\hat{Z}_{t}^{\e,\xi}})\|^{2}\bigg]\\ 
 &&+\frac{1}{\e}
 \mE\bigg[2\<Z_{t}^{\e,\xi}-\hat{Z}_{t}^{\e,\xi}, b_{2}(X_{t}^{\e},\sL^{\mP}_{X_{t}^{\e}},\hat{Z}_{t}^{\e,\xi},\sL^{\mP}_{\hat{Z}_{t}^{\e,\xi}})
 -b_{2}(X_{t(\d)}^{\e},\sL^{\mP}_{X_{t(\d)}^{\e}},\hat{Z}_{t}^{\e,\xi},\sL^{\mP}_{\hat{Z}_{t}^{\e,\xi}})\>\\
 &&+(2p+1)\|\s_{2}(X_{t}^{\e},\sL^{\mP}_{X_{t}^{\e}},\hat{Z}_{t}^{\e,\xi},\sL^{\mP}_{\hat{Z}_{t}^{\e,\xi}})
 -\s_{2}(X_{t(\d)}^{\e},\sL^{\mP}_{X_{t(\d)}^{\e}},\hat{Z}_{t}^{\e,\xi},\sL^{\mP}_{\hat{Z}_{t}^{\e,\xi}})\|^{2}\bigg] \\
 &\leq&\frac{1}{\e}\mE\Big(-(\b_{1}-L_{b_{2},\s_2})
 |Z_{t}^{\e,\xi}-\hat{Z}_{t}^{\e,\xi}|^{2}+\b_{2}\mW_2^{2}(\sL^{\mP}_{Z_{t}^{\e,\xi}},\sL^{\mP}_{\hat{Z}_{t}^{\e,\xi}})\Big)\\
 &&+\frac{C}{\e}\mE
 \(|X_{t}^{\e}-X_{t(\d)}^{\e}|^{2}+\mW_2^{2}(\sL^{\mP}_{X_{t}^{\e}},\sL^{\mP}_{X_{t(\d)}^{\e}})\)\\
 &\leq&\frac{-(\b_1-\b_2-L_{b_{2},\s_2})}{\e}\mE|Z_{t}^{\e,\xi}-\hat{Z}_{t}^{\e,\xi}|^{2}+\frac{C}{\e}\mE|X_{t}^{\e}-X_{t(\d)}^{\e}|^{2}.
\de
Then it follows from the comparison theorem and (\ref{barx}) that
\ce
\mE|Z_{t}^{\e,\xi}-\hat{Z}_{t}^{\e,\xi}|^{2}\leq \frac{C}{\e}({\d}^{2}+{\d})\int_0^t e^{\frac{-(\b_1-\b_2-L_{b_{2},\s_2})}{\e}(t-s)}\dif s
\leq\frac{C}{\b_1-\b_2-L_{b_{2},\s_2}}({\d}^{2}+{\d}).
\de

Next, by the similar deduction to the above, we have that
\ce
\mE|Z_{t}^{\e,z_0,\sL^{\mP}_\xi}-\hat{Z}_{t}^{\e,z_0,\sL^{\mP}_\xi}|^{2}\leq \frac{C}{\b_1-L_{b_{2},\s_2}}({\d}^{2}+{\d}).
\de
The proof is complete.
\end{proof}

Finally, we estimate $\mE\sup\limits_{t\in[0,T]}|X_t^\e-\hat{X}_t^\e|^2$. By (\ref{Eq1}) (\ref{hatx}), it holds that
\be
&&\mE\sup\limits_{t\in[0,T]}|X_t^\e-\hat{X}_t^\e|^2\no\\
&\leq& \mE\sup\limits_{t\in[0,T]}\left|\int_0^t\left[b_1(X^\e_{s},\sL^{\mP}_{X^\e_{s}},Z_{s}^{\e,z_0,\sL^{\mP}_\xi},\sL^{\mP}_{Z_{s}^{\e,\xi}})-b_1(X^\e_{s(\d)},\sL^{\mP}_{X^\e_{s(\d)}},\hat{Z}_{s}^{\e,z_0,\sL^{\mP}_\xi},\sL^{\mP}_{\hat{Z}_{s}^{\e,\xi}})\right]\dif s\right|^2\no\\
&\leq&TL_{b_1,\s_1}\mE\int_0^T(|X^\e_{s}-X^\e_{s(\d)}|^2+\mW^2_2(\sL^{\mP}_{X^\e_{s}},\sL^{\mP}_{X^\e_{s(\d)}})+|Z_{s}^{\e,z_0,\sL^{\mP}_\xi}-\hat{Z}_{s}^{\e,z_0,\sL^{\mP}_\xi}|^2\no\\
&&\qquad +\mW^2_2(\sL^{\mP}_{Z_{s}^{\e,\xi}},\sL^{\mP}_{\hat{Z}_{s}^{\e,\xi}}))\dif s\no\\
&\leq&C({\d}^{2}+{\d}),
\label{xehatxe}
\ee
where we use (\ref{barx}) and Lemma \ref{zehatze}.

\subsection{Some estimates for the frozen equation (\ref{Eq2})}\label{froequ}

Consider the following classical SDE related with the frozen equation (\ref{Eq2}):
\be\left\{\begin{array}{l}
\dif Z_{t}^{x,\mu,z_0,\sL^{\mP}_\xi}=b_{2}(x,\mu,Z_{t}^{x,\mu,z_0,\sL^{\mP}_\xi},\sL^{\mP}_{Z_{t}^{x,\mu,\xi}})\dif t+\s_{2}(x,\mu,Z_{t}^{x,\mu,z_0,\sL^{\mP}_\xi},\sL^{\mP}_{Z_{t}^{x,\mu,\xi}})\dif W_{t},\\
Z_{0}^{x,\mu,z_0,\sL^{\mP}_\xi}=z_0, \quad 0 \leq t \leq T.
\end{array}
\right.
\label{Eq21}
\ee
Under $(\mathbf{H}^1_{b_{2}, \s_{2}})$, we know that the above equation has a unique strong solution $Z_{\cdot}^{x,\mu,z_0,\sL^{\mP}_\xi}$. Moreover, about $Z_{\cdot}^{x,\mu,z_0,\sL^{\mP}_\xi}$, we have the following estimates.

\bl\label{xmuzxi}
Under these assumptions $(\mathbf{H}^1_{b_{2}, \s_{2}})$-$(\mathbf{H}^2_{b_{2}, \s_{2}})$, it holds that for any $t\in[0,T], x\in\mR^n, \mu\in\cP_2(\mR^n)$
\be
\mE|Z_{t}^{x, \mu,z_0,\sL^{\mP}_\xi}|^{2}\leq|z_0|^2e^{-\a_1t}+C(\|\sL^{\mP}_\xi\|^{2}+1+|x|^{2}+\|\mu\|^{2}).
\label{memu2}
\ee
\el
\begin{proof}
First of all, we estimate $Z_{\cdot}^{x, \mu,\xi}$. Applying the It\^{o} formula to $|Z_{t}^{x, \mu,\xi}|^{2}$ and taking the expectation, we get that
\ce
\mE|Z_{t}^{x, \mu, \xi}|^{2}
&=&\mE|\xi|^{2}+2\mE\int_{0}^{t}\<Z_{s}^{x, \mu, \xi}, b_{2}(x,\mu, Z_{s}^{x, \mu, \xi}, \sL^{\mP}_{Z_{s}^{x, \mu, \xi}})\>\dif s\\
&&+\mE\int_{0}^{t}\|\s_{2}(x, \mu, Z_{s}^{x, \mu, \xi}, \sL^{\mP}_{Z_{s}^{x, \mu, \xi}})\|^{2}\dif s,
\de
and 
\ce
\frac{\dif}{\dif t}\mE|Z_{t}^{x, \mu, \xi}|^{2}&=&\mE\(2\<Z_{t}^{x, \mu, \xi}, b_{2}(x,\mu, Z_{t}^{x, \mu, \xi}, \sL^{\mP}_{Z_{t}^{x, \mu, \xi}})\>+\|\s_{2}(x, \mu, Z_{t}^{x, \mu, \xi}, \sL^{\mP}_{Z_{t}^{x, \mu, \xi}})\|^2\)\\
&\leq&\mE\(-\a_{1}|Z_{t}^{x, \mu, \xi}|^{2}+\a_{2}\|\sL^{\mP}_{Z_{t}^{x, \mu, \xi}}\|^{2}+C(1+|x|^{2}+\|\mu\|^{2})\)\\
&=&-(\a_1-\a_2)\mE|Z_{t}^{x, \mu, \xi}|^{2}+C(1+|x|^{2}+\|\mu\|^{2}).
\de
By the comparison theorem, it holds that
\be
\mE|Z_{t}^{x, \mu, \xi}|^{2}&\leq&\mE|\xi|^{2}e^{-(\a_1-\a_2) t}+C(1+|x|^{2}+\|\mu\|^{2})\int_{0}^{t}e^{-(\a_1-\a_2)(t-s)}\dif s\no\\
&\leq&\mE|\xi|^{2}e^{-(\a_1-\a_2) t}+C(1+|x|^{2}+\|\mu\|^{2}).
\label{zxmuxi}
\ee

Next, we deal with $Z_{\cdot}^{x,\mu,z_0,\sL^{\mP}_\xi}$. By the same deduction to the above, it holds that
\ce
\mE|Z_{t}^{x, \mu,z_0,\sL^{\mP}_\xi}|^{2}\leq |z_0|^2e^{-\a_1t}+C(\|\sL^{\mP}_\xi\|^{2}+1+|x|^{2}+\|\mu\|^{2}).
\de
The proof is complete.
\end{proof}

\bl\label{xmuzxi12}
Suppose that $(\mathbf{H}^{1}_{b_{2}, \s_{2}})$ $(\mathbf{H}^2_{b_{2}, \s_{2}})$ hold. Then it holds that for any $x_i\in\mR^n, \mu_i\in\cP_2(\mR^n), z_i\in\mR^m, \zeta_i\in L^{2}(\Omega,\sF_0,\mP;\mR^m), i=1,2$,
\ce
&&\mE|Z_{t}^{x_1,\mu_1,z_1,\sL^{\mP}_{\zeta_1}}-Z_{t}^{x_2,\mu_2,z_2,\sL^{\mP}_{\zeta_2}}|^{2}\\
&\leq& |z_1-z_2|^2e^{-(\b_{1}-L_{b_{2},\s_2})t}+\mE|\zeta_{1}-\zeta_{2}|^{2}\left(e^{\b_2 t}-1\right)e^{-(\b_{1}-L_{b_{2},\s_2})t}\\
&&+C(|x_1-x_2|^2+\mW^2_2(\mu_1,\mu_2))\frac{1-e^{-(\b_{1}-L_{b_{2},\s_2})t}}{\b_{1}-L_{b_{2},\s_2}}.
\de
\el
\begin{proof}
First of all, we compute $\mE|Z_{t}^{x_1,\mu_1,\zeta_{1}}-Z_{t}^{x_2,\mu_2,\zeta_{2}}|^{2}$. Note that $Z_{t}^{x_1,\mu_1,\zeta_{1}}$ and $Z_{t}^{x_2,\mu_2,\zeta_{2}}$ solve Eq.(\ref{Eq2}) with initial values $\zeta_{1}$ and $\zeta_{2}$, respectively, i.e.
\ce
&&Z_{t}^{x_1,\mu_1,\zeta_{1}}-Z_{t}^{x_2,\mu_2,\zeta_{2}}\\
&=&\zeta_{1}-\zeta_{2}+\int_{0}^{t}\(b_{2}(x_1,\mu_1,Z_{s}^{x_1,\mu_1,\zeta_{1}},\sL^{\mP}_{Z_{s}^{x_1,\mu_1,\zeta_{1}}})
-b_{2}(x_2,\mu_2,Z_{s}^{x_2,\mu_2,\zeta_{2}},\sL^{\mP}_{Z_{s}^{x_2,\mu_2,\zeta_{2}}})\)\dif s\\
&&+\int_{0}^{t}\(\s_{2}(x_1,\mu_1,Z_{s}^{x_1,\mu_1,\zeta_{1}},\sL^{\mP}_{Z_{s}^{x_1,\mu_1,\zeta_{1}}})
-\s_{2}(x_2,\mu_2,Z_{s}^{x_2,\mu_2,\zeta_{2}},\sL^{\mP}_{Z_{s}^{x_2,\mu_2,\zeta_{2}}})\)\dif W_{s}.
\de
Applying the It\^{o} formula to $|Z_{t}^{x_1,\mu_1,\zeta_{1}}-Z_{t}^{x_2,\mu_2,\zeta_{2}}|^{2}$ and taking expectation on two sides, we obtain that
\ce
&&\mE|Z_{t}^{x_1,\mu_1,\zeta_{1}}-Z_{t}^{x_2,\mu_2,\zeta_{2}}|^{2}\\
&=&\mE|\zeta_{1}-\zeta_{2}|^{2}+2\mE\int_{0}^{t}\< Z_{s}^{x_1,\mu_1,\zeta_{1}}-Z_{s}^{x_2,\mu_2,\zeta_{2}}, b_{2}(x_1,\mu_1,Z_{s}^{x_1,\mu_1,\zeta_{1}},\sL^{\mP}_{Z_{s}^{x_1,\mu_1,\zeta_{1}}})\\
&&\quad\quad-b_{2}(x_2,\mu_2,Z_{s}^{x_2,\mu_2,\zeta_{2}},\sL^{\mP}_{Z_{s}^{x_2,\mu_2,\zeta_{2}}})\>
\dif s\\
&&+\mE\int_{0}^{t}\|\s_{2}(x_1,\mu_1,Z_{s}^{x_1,\mu_1,\zeta_{1}},\sL^{\mP}_{Z_{s}^{x_1,\mu_1,\zeta_{1}}})
-\s_{2}(x_2,\mu_2,Z_{s}^{x_2,\mu_2,\zeta_{2}},\sL^{\mP}_{Z_{s}^{x_2,\mu_2,\zeta_{2}}})\|^{2}\dif s,
\de
and
\ce
&&\frac{\dif}{\dif t}\mE|Z_{t}^{x_1,\mu_1,\zeta_{1}}-Z_{t}^{x_2,\mu_2,\zeta_{2}}|^{2}\\
&=&\mE\(2\< Z_{t}^{x_1,\mu_1,\zeta_{1}}-Z_{t}^{x_2,\mu_2,\zeta_{2}}, b_{2}(x_1,\mu_1,Z_{t}^{x_1,\mu_1,\zeta_{1}},\sL^{\mP}_{Z_{t}^{x_1,\mu_1,\zeta_{1}}})-b_{2}(x_2,\mu_2,Z_{t}^{x_2,\mu_2,\zeta_{2}},\sL^{\mP}_{Z_{t}^{x_2,\mu_2,\zeta_{2}}})\>\)\\
&&+\mE\|\s_{2}(x_1,\mu_1,Z_{t}^{x_1,\mu_1,\zeta_{1}},\sL^{\mP}_{Z_{t}^{x_1,\mu_1,\zeta_{1}}})
-\s_{2}(x_2,\mu_2,Z_{t}^{x_2,\mu_2,\zeta_{2}},\sL^{\mP}_{Z_{t}^{x_2,\mu_2,\zeta_{2}}})\|^{2}\\
&\leq&\mE\(2\< Z_{t}^{x_1,\mu_1,\zeta_{1}}-Z_{t}^{x_2,\mu_2,\zeta_{2}}, b_{2}(x_1,\mu_1,Z_{t}^{x_1,\mu_1,\zeta_{1}},\sL^{\mP}_{Z_{t}^{x_1,\mu_1,\zeta_{1}}})-b_{2}(x_1,\mu_1,Z_{t}^{x_2,\mu_2,\zeta_{2}},\sL^{\mP}_{Z_{t}^{x_2,\mu_2,\zeta_{2}}})\>\)\\
&&+(2p+1)\mE\|\s_{2}(x_1,\mu_1,Z_{t}^{x_1,\mu_1,\zeta_{1}},\sL^{\mP}_{Z_{t}^{x_1,\mu_1,\zeta_{1}}})
-\s_{2}(x_1,\mu_1,Z_{t}^{x_2,\mu_2,\zeta_{2}},\sL^{\mP}_{Z_{t}^{x_2,\mu_2,\zeta_{2}}})\|^{2}\\
&&+\mE\(2\< Z_{t}^{x_1,\mu_1,\zeta_{1}}-Z_{t}^{x_2,\mu_2,\zeta_{2}}, b_{2}(x_1,\mu_1,Z_{t}^{x_2,\mu_2,\zeta_{2}},\sL^{\mP}_{Z_{t}^{x_2,\mu_2,\zeta_{2}}})-b_{2}(x_2,\mu_2,Z_{t}^{x_2,\mu_2,\zeta_{2}},\sL^{\mP}_{Z_{t}^{x_2,\mu_2,\zeta_{2}}})\>\)\\
&&+(2p+1)\mE\|\s_{2}(x_1,\mu_1,Z_{t}^{x_2,\mu_2,\zeta_{2}},\sL^{\mP}_{Z_{t}^{x_2,\mu_2,\zeta_{2}}})
-\s_{2}(x_2,\mu_2,Z_{t}^{x_2,\mu_2,\zeta_{2}},\sL^{\mP}_{Z_{t}^{x_2,\mu_2,\zeta_{2}}})\|^{2}\\
&\leq&\mE\(-(\b_{1}-L_{b_{2},\s_2})|Z_{t}^{x_1,\mu_1,\zeta_{1}}-Z_{t}^{x_2,\mu_2,\zeta_{2}}|^{2}
+\b_{2}\mW_2^{2}(\sL^{\mP}_{Z_{t}^{x_1,\mu_1,\zeta_{1}}}, \sL^{\mP}_{Z_{t}^{x_2,\mu_2,\zeta_{2}}})\)\\
&&+C(|x_1-x_2|^2+\mW^2_2(\mu_1,\mu_2))\\
&\leq&-(\b_{1}-\b_{2}-L_{b_{2},\s_2})\mE|Z_{t}^{x_1,\mu_1,\zeta_{1}}-Z_{t}^{x_2,\mu_2,\zeta_{2}}|^{2}+C(|x_1-x_2|^2+\mW^2_2(\mu_1,\mu_2)).
\de
By the comparison theorem, it holds that
$$
\mE|Z_{t}^{x_1,\mu_1,\zeta_{1}}-Z_{t}^{x_2,\mu_2,\zeta_{2}}|^{2}\leq\mE|\zeta_{1}-\zeta_{2}|^{2}e^{-(\b_{1}-\b_{2}-L_{b_{2},\s_2})t}+C(|x_1-x_2|^2+\mW^2_2(\mu_1,\mu_2)).
$$

Next, we investigate $\mE|Z_{t}^{x_1,\mu_1,z_1,\sL^{\mP}_{\zeta_1}}-Z_{t}^{x_2,\mu_2,z_2,\sL^{\mP}_{\zeta_2}}|^{2}$. The It\^o formula yields that
\ce
&&\mE|Z_{t}^{x_1,\mu_1,z_1,\sL^{\mP}_{\zeta_1}}-Z_{t}^{x_2,\mu_2,z_2,\sL^{\mP}_{\zeta_2}}|^{2}e^{\lambda t}\\
&=&|z_1-z_2|^{2}+\l \mE\int_{0}^{t}e^{\l s}|Z_{s}^{x_1,\mu_1,z_1,\sL^{\mP}_{\zeta_1}}-Z_{s}^{x_2,\mu_2,z_2,\sL^{\mP}_{\zeta_2}}|^{2}\dif s\\
&&+2\mE\int_{0}^{t}e^{\l s}\< Z_{s}^{x_1,\mu_1,z_1,\sL^{\mP}_{\zeta_1}}-Z_{s}^{x_2,\mu_2,z_2,\sL^{\mP}_{\zeta_2}}, b_{2}(x_1,\mu_1,Z_{s}^{x_1,\mu_1,z_1,\sL^{\mP}_{\zeta_1}},\sL^{\mP}_{Z_{s}^{x_1,\mu_1,\zeta_{1}}})\\
&&\quad\quad-b_{2}(x_2,\mu_2,Z_{s}^{x_2,\mu_2,z_2,\sL^{\mP}_{\zeta_2}},\sL^{\mP}_{Z_{s}^{x_2,\mu_2,\zeta_{2}}})\>
\dif s\\
&&+\mE\int_{0}^{t}e^{\l s}\|\s_{2}(x_1,\mu_1,Z_{s}^{x_1,\mu_1,z_1,\sL^{\mP}_{\zeta_1}},\sL^{\mP}_{Z_{s}^{x_1,\mu_1,\zeta_{1}}})
-\s_{2}(x_2,\mu_2,Z_{s}^{x_2,\mu_2,z_2,\sL^{\mP}_{\zeta_2}},\sL^{\mP}_{Z_{s}^{x_2,\mu_2,\zeta_{2}}})\|^{2}\dif s\\
&\leq&|z_1-z_2|^{2}+\l \mE\int_{0}^{t}e^{\l s}|Z_{s}^{x_1,\mu_1,z_1,\sL^{\mP}_{\zeta_1}}-Z_{s}^{x_2,\mu_2,z_2,\sL^{\mP}_{\zeta_2}}|^{2}\dif s\\
&&-(\b_{1}-L_{b_{2},\s_2})\mE\int_{0}^{t}e^{\l s}|Z_{s}^{x_1,\mu_1,z_1,\sL^{\mP}_{\zeta_1}}-Z_{s}^{x_2,\mu_2,z_2,\sL^{\mP}_{\zeta_2}}|^{2}\dif s\\
&&+\b_{2}\int_{0}^{t}e^{\l s}\mW_2^{2}(\sL^{\mP}_{Z_{s}^{x_1,\mu_1,\zeta_{1}}}, \sL^{\mP}_{Z_{s}^{x_2,\mu_2,\zeta_{2}}})\)\dif s+C(|x_1-x_2|^2+\mW^2_2(\mu_1,\mu_2))\int_{0}^{t}e^{\l s}\dif s\\
&\leq&|z_1-z_2|^{2}+\b_{2}\mE|\zeta_{1}-\zeta_{2}|^{2}\int_{0}^{t}e^{\l s}e^{-(\b_{1}-\b_{2}-L_{b_{2},\s_2})s}\dif s+C(|x_1-x_2|^2+\mW^2_2(\mu_1,\mu_2))\int_{0}^{t}e^{\l s}\dif s,
\de
where $\l:=\b_{1}-L_{b_{2},\s_2}$. Then simple calculation implies that
\ce
&&\mE|Z_{t}^{x_1,\mu_1,z_1,\sL^{\mP}_{\zeta_1}}-Z_{t}^{x_2,\mu_2,z_2,\sL^{\mP}_{\zeta_2}}|^{2}\\
&\leq& |z_1-z_2|^2e^{-(\b_{1}-L_{b_{2},\s_2})t}+\mE|\zeta_{1}-\zeta_{2}|^{2}\left(e^{\b_2 t}-1\right)e^{-(\b_{1}-L_{b_{2},\s_2})t}\\
&&+C(|x_1-x_2|^2+\mW^2_2(\mu_1,\mu_2))\frac{1-e^{-(\b_{1}-L_{b_{2},\s_2})t}}{\b_{1}-L_{b_{2},\s_2}},
\de
which completes the proof.
\end{proof}

\bl\label{emb1}
 Suppose that $(\mathbf{H}^1_{b_{1}, \s_{1}})$ $(\mathbf{H}^{1}_{b_{2}, \s_{2}})$ $(\mathbf{H}^2_{b_{2}, \s_{2}})$ hold. Then there exists a constant $C>0$ such that for any $t\in[0, T]$, $x\in\mR^n, \mu\in\cP_2(\mR^n)$
\be
|\mE b_{1}(x, \mu, Z_{t}^{x,\mu,z_0,\sL^{\mP}_\xi},\sL_{Z_{t}^{x,\mu,\xi}})-\bar{b}_{1}(x, \mu)|^2\leq Ce^{-(\b_{1}-L_{b_{2},\s_2})t}(\|\sL^{\mP}_\xi\|+1+|x|+\|\mu\|+|z_0|)^2.
\label{meu2}
\ee
\el
\begin{proof}
Based on $(\mathbf{H}^1_{b_{1}, \s_{1}})$ and Lemma \ref{xmuzxi12}, one can obtain that
\ce
&&|\mE b_{1}(x, \mu, Z_{t}^{x,\mu,z_0,\sL^{\mP}_\xi},\sL^{\mP}_{Z_{t}^{x,\mu,\xi}})-\bar{b}_{1}(x, \mu)|^2\\
&=&\Big|\mE b_{1}(x, \mu, Z_{t}^{x,\mu,z_0,\sL^{\mP}_\xi},\sL^{\mP}_{Z_{t}^{x,\mu,\xi}})-\int_{\mR^{m}\times\cP_2(\mR^{m})}b_{1}(x,\mu,y,\nu)\eta^{x,\mu}\times\d_{\eta^{x,\mu}}(\dif y,\dif \nu)\Big|^2\\
&\overset{\sL^{\mP}_{\zeta}=\nu}{=}&\Big|\mE b_{1}(x, \mu, Z_{t}^{x,\mu,z_0,\sL^{\mP}_\xi},\sL^{\mP}_{Z_{t}^{x,\mu,\xi}})-\int_{\mR^{m}\times\cP_2(\mR^{m})}\mE b_{1}(x, \mu, Z_{t}^{x,\mu,y,\nu},\sL^{\mP}_{Z_{t}^{x,\mu,\zeta}})\eta^{x,\mu}\times\d_{\eta^{x,\mu}}(\dif y,\dif \nu)\Big|^2\\
&\leq&\int_{\mR^{m}\times\cP_2(\mR^{m})}\mE |b_{1}(x, \mu, Z_{t}^{x,\mu,z_0,\sL^{\mP}_\xi},\sL^{\mP}_{Z_{t}^{x,\mu,\xi}})- b_{1}(x, \mu, Z_{t}^{x,\mu,y,\nu},\sL^{\mP}_{Z_{t}^{x,\mu,\zeta}})|^2\eta^{x,\mu}\times\d_{\eta^{x,\mu}}(\dif y,\dif \nu)\\
&\overset{\sL^{\mP}_{\zeta^{x,\mu}}=\eta^{x,\mu}}{\leq}& L_{b_{1},\s_{1}}\int_{\mR^{m}}\(\mE |Z_{t}^{x,\mu,z_0,\sL^{\mP}_\xi}-Z_{t}^{x,\mu,y,\eta^{x,\mu}}|^2+\mW^2_2(\sL^{\mP}_{Z_{t}^{x,\mu,\xi}},\sL^{\mP}_{Z_{t}^{x,\mu,\zeta^{x,\mu}}})\)\eta^{x,\mu}(\dif y)\\
&\leq& C\int_{\mR^{m}}|z_0-y|^2e^{-(\b_{1}-L_{b_{2},\s_2})t}\eta^{x,\mu}(\dif y)+C\mW^2_2(\sL^{\mP}_{\xi},\eta^{x,\mu})\left(e^{\b_2 t}-1\right)e^{-(\b_{1}-L_{b_{2},\s_2})t}\\
&&+C\mE|\xi-\zeta^{x,\mu}|^2 e^{-(\b_{1}-\b_{2}-L_{b_{2},\s_2})t}\\
&\overset{(\ref{zxmuxi})}{\leq}&Ce^{-(\b_{1}-L_{b_{2},\s_2})t}(\|\sL^{\mP}_{\xi}\|+1+|x|+\|\mu\|+|z_0|)^2.
\de
The proof is complete. 
\end{proof}

\subsection{Some estimates for the average equation (\ref{Eq3})}\label{aveequ}

\bl
Suppose that $(\mathbf{H}^{1}_{b_{1}, \s_{1}})$ $(\mathbf{H}^{1}_{b_{2}, \s_{2}})$-$(\mathbf{H}^{2}_{b_{2}, \s_{2}})$ hold. Then Eq.(\ref{Eq3}) has a unique strong solution $\bar{X_{\cdot}}$. Moreover,
\be
\mE\left(\sup\limits_{t\in[0,T]}|\bar{X}_{t}|^{2p+2}\right)\leq C(1+\mE|\varrho|^{2p+2}).
\label{barxb}
\ee
\el
\begin{proof}
First of all, we justify that for any $x_i\in\mR^n, \mu_i\in\cP_2(\mR^n)$, $i=1,2$
\ce
|\bar{b}_{1}(x_{1},\mu_{1})-\bar{b}_{1}(x_{2},\mu_{2})|^2\leq L_{\bar{b}_{1}}\(|x_{1}-x_{2}|^{2}+\mW_2^{2}(\mu_{1}, \mu_{2})\),
\de
where $L_{\bar{b}_{1}}>0$ is a constant. Indeed, by Lemma \ref{xmuzxi12}, Lemma \ref{emb1} and $(\mathbf{H}^{1}_{b_{1}, \s_{1}})$, it holds that
\ce
&&|\bar{b}_{1}(x_{1},\mu_{1})-\bar{b}_{1}(x_{2},\mu_{2})|^2\\
&\leq& 3|\bar{b}_{1}(x_{1},\mu_{1})-\mE b_1(x_{1},\mu_{1},Z_{t}^{x_1,\mu_1,z_0,\sL^{\mP}_{\xi}},\sL^{\mP}_{Z_{t}^{x_1,\mu_1,\xi}})|^2\\
&&+3|\bar{b}_{1}(x_{2},\mu_{2})-\mE b_1(x_{2},\mu_{2},Z_{t}^{x_2,\mu_2,z_0,\sL^{\mP}_{\xi}},\sL^{\mP}_{Z_{t}^{x_2,\mu_2,\xi}})|^2\\
&&+3|\mE b_1(x_{1},\mu_{1},Z_{t}^{x_1,\mu_1,z_0,\sL^{\mP}_{\xi}},\sL^{\mP}_{Z_{t}^{x_1,\mu_1,\xi}})-\mE b_1(x_{2},\mu_{2},Z_{t}^{x_2,\mu_2,z_0,\sL^{\mP}_{\xi}},\sL^{\mP}_{Z_{t}^{x_2,\mu_2,\xi}})|^2\\
&\leq&Ce^{-(\b_{1}-L_{b_{2},\s_2})t}(1+|x_1|^2+|x_2|^2+\|\mu_1\|^2+\|\mu_2\|^2+|z_0|^2+\|\sL^{\mP}_{\xi}\|^2)\\
&&+C\(|x_{1}-x_{2}|^{2}+\mW_2^{2}(\mu_{1},\mu_{2})),
\de
which implies the required result as $t\rightarrow \infty$. Thus, from \cite[Theorem 2.1]{WangFY} it follows that Eq.(\ref{Eq3}) has a unique strong solution $\bar{X_{\cdot}}$. 

Next, by similar deduction to that for $X^\e$ in Lemma \ref{xtztc}, we have (\ref{barxb}). The proof is complete.
\end{proof}

\bl\label{2orde}
Suppose that assumptions $(\mathbf{H}^{1}_{b_{1}, \s_{1}})$ $(\mathbf{H}^{1}_{b_{2}, \s_{2}})$-$(\mathbf{H}^{2}_{b_{2}, \s_{2}})$ hold. Then there exists a constant $C>0$ such that 
\ce
\mE\(\sup_{0\leq t\leq T}|\hat{X}_{t}^{\e}-\bar{X}_{t}|^{2}\)\leq C\(\frac{\e}{\d}+{\d}^{2}+{\d}\).
\de
\el
\begin{proof}
{\bf Step 1.} We estimate $\hat{X}_{t}^{\e}-\bar{X}_{t}$.

Note that
\ce
\hat{X}_{t}^{\e}-\bar{X}_{t}
&=&\int_{0}^{t}\(b_{1}(X_{s(\d)}^{\e},\sL^{\mP}_{X_{s(\d)}^{\e}},\hat{Z}_{s}^{\e,z_0,\sL^{\mP}_{\xi}},\sL^{\mP}_{\hat{Z}_{s}^{\e,\xi}})-\bar{b}_{1}(\bar{X}_{s},\sL^{\mP}_{\bar{X}_{s}})\)\dif s\\
&&+\int_{0}^{t}\(\s_{1}(X_{s}^{\e},\sL^{\mP}_{X_{s}^{\e}})-\s_{1}(\bar{X}_{s},\sL^{\mP}_{\bar{X}_{s}})\)\dif B_{s}.
\de
Thus, based on the BDG inequality, we get that
\ce
&&\mE\(\sup_{0\leq t\leq T}|\hat{X}_{t}^{\e}-\bar{X}_{t}|^{2}\)\\
&\leq& 2\Bigg(\mE\sup_{0\leq t\leq T}
\Big|\int_{0}^{t}\(b_{1}(X_{s(\d)}^{\e},\sL^{\mP}_{X_{s(\d)}^{\e}},\hat{Z}_{s}^{\e,z_0,\sL^{\mP}_{\xi}},\sL^{\mP}_{\hat{Z}_{s}^{\e,\xi}})-\bar{b}_{1}(\bar{X}_{s},\sL^{\mP}_{\bar{X}_{s}})\)
\dif s\Big|^{2}\Bigg)\\
&&+2\Bigg(\mE\sup_{0\leq t\leq T}
\Big|\int_{0}^{t}\(\s_{1}(X_{s}^{\e},\sL^{\mP}_{X_{s}^{\e}})-\s_{1}(\bar{X}_{s},\sL^{\mP}_{\bar{X}_{s}})\)\dif B_{s}\Big|^{2}\Bigg)\\
 &\leq&6\Bigg(\mE\sup_{0\leq t\leq T}
 \Big|\int_{0}^{t}\(b_{1}(X_{s(\d)}^{\e},\sL^{\mP}_{X_{s(\d)}^{\e}},\hat{Z}_{s}^{\e,z_0,\sL^{\mP}_{\xi}},\sL^{\mP}_{\hat{Z}_{s}^{\e,\xi}})-\bar{b}_{1}(X_{s(\d)}^{\e},\sL^{\mP}_{X_{s(\d)}^{\e}})\)
 \dif s\Big|^{2}\Bigg)\\
 &&+6\Bigg(\mE\sup_{0\leq t\leq T}
 \Big|\int_{0}^{t}\(\bar{b}_{1}(X_{s(\d)}^{\e},\sL^{\mP}_{X_{s(\d)}^{\e}})-\bar{b}_{1}(X_{s}^{\e},\sL^{\mP}_{X_{s}^{\e}})\)
 \dif s\Big|^{2}\Bigg)\\
 &&+6\Bigg(\mE\sup_{0\leq t\leq T}
 \Big|\int_{0}^{t}\(\bar{b}_{1}(X_{s}^{\e},\sL^{\mP}_{X_{s}^{\e}})-\bar{b}_{1}(\bar{X}_{s},\sL^{\mP}_{\bar{X}_{s}})\)\dif s\Big|^{2})\\
 &&+8
 \int_{0}^{T}\mE\Big\|\s_{1}(X_{s}^{\e},\sL^{\mP}_{X_{s}^{\e}})-\s_{1}(\bar{X}_{s},\sL^{\mP}_{\bar{X}_{s}})\Big\|^{2}\dif s.
\de
Then from the H\"{o}lder inequality and $(\mathbf{H}^1_{b_{1}, \s_{1}})$, it follows that
\be
&&\mE\(\sup_{0\leq t\leq T}|\hat{X}_{t}^{\e}-\bar{X}_{t}|^{2}\)\no\\
&\leq&6\Bigg(\mE\sup_{0\leq t\leq T}
\Big|\int_{0}^{t}\(b_{1}(X_{s(\d)}^{\e},\sL^{\mP}_{X_{s(\d)}^{\e}},\hat{Z}_{s}^{\e,z_0,\sL^{\mP}_{\xi}},\sL^{\mP}_{\hat{Z}_{s}^{\e,\xi}})-\bar{b}_{1}(X_{s(\d)}^{\e},\sL^{\mP}_{X_{s(\d)}^{\e}})\)
\dif s\Big|^{2}\Bigg)\no\\
&&+6T\Bigg(\mE\sup_{0\leq t\leq T}
\int_{0}^{t}L_{\bar{b}_{1}}\Big(|X_{s}^{\e}-X_{s(\d)}^{\e}|^{2}+\mW_2^{2}(\sL^{\mP}_{X_{s}^{\e}},\sL^{\mP}_{X_{s(\d)}^{\e}})\Big)
\dif s\Bigg)\no\\
&&+6T\Bigg(\mE\sup_{0\leq t\leq T}
\int_{0}^{t}L_{\bar{b}_{1}}\Big(|X_{s}^{\e}-\bar{X}_{s}|^{2}+\mW_2^{2}(\sL^{\mP}_{X_{s}^{\e}},\sL^{\mP}_{\bar{X}_{s}})\Big)
\dif s\Bigg)\no\\
&&+8
\int_{0}^{T}\mE L_{b_{1},\s_{1}}\Big(|X_{s}^{\e}-\bar{X}_{s}|^{2}+\mW_2^{2}(\sL^{\mP}_{X_{s}^{\e}},\sL^{\mP}_{\bar{X}_{s}})\Big)
\dif s\no\\
&\leq&6\Bigg(\mE\sup_{0\leq t\leq T}
\Big|\int_{0}^{t}\(b_{1}(X_{s(\d)}^{\e},\sL^{\mP}_{X_{s(\d)}^{\e}},\hat{Z}_{s}^{\e,z_0,\sL^{\mP}_{\xi}},\sL^{\mP}_{\hat{Z}_{s}^{\e,\xi}})
-\bar{b}_{1}(X_{s(\d)}^{\e},\sL^{\mP}_{X_{s(\d)}^{\e}})\)\dif s\Big|^{2}\Bigg)\no\\
 &&+12TL_{\bar{b}_1}\int_{0}^{T}\mE|X_{s}^{\e}-X_{s(\d)}^{\e}|^{2}\dif s
+(12TL_{\bar{b}_1}+16L_{b_{1},\s_{1}})\int_{0}^{T}\mE|X_{s}^{\e}-\bar{X}_{s}|^{2}\dif s\no\\
&\leq&6\Bigg(\mE\sup_{0\leq t\leq T}
\Big|\int_{0}^{t}\(b_{1}(X_{s(\d)}^{\e},\sL^{\mP}_{X_{s(\d)}^{\e}},\hat{Z}_{s}^{\e,z_0,\sL^{\mP}_{\xi}},\sL^{\mP}_{\hat{Z}_{s}^{\e,\xi}})
-\bar{b}_{1}(X_{s(\d)}^{\e},\sL^{\mP}_{X_{s(\d)}^{\e}})\)\dif s\Big|^{2}\Bigg)\no\\
 &&+12TL_{\bar{b}_1}\int_{0}^{T}\mE|X_{s}^{\e}-X_{s(\d)}^{\e}|^{2}\dif s
+2(12TL_{\bar{b}_1}+16L_{b_{1},\s_{1}})\int_{0}^{T}\mE|X_{s}^{\e}-\hat{X}^\e_{s}|^{2}\dif s\no\\
&&+2(12TL_{\bar{b}_1}+16L_{b_{1},\s_{1}})\int_{0}^{T}\mE|\hat{X}^\e_{s}-\bar{X}_{s}|^{2}\dif s\no\\
&=:&I_1+I_2+I_3+2(12TL_{\bar{b}_1}+16L_{b_{1},\s_{1}})\int_{0}^{T}\mE|\hat{X}^\e_{s}-\bar{X}_{s}|^{2}\dif s.
\label{meb4}
\ee

By the deduction in {\bf Step 2}, we know that 
\be
I_1\leq C(\frac{\e}{\d}+\d).
\label{b4le}
\ee
And (\ref{barx}) (\ref{xehatxe}) imply that
\be
I_2+I_3\leq C(\d^2+\d).
\label{i2i3}
\ee

Thus, inserting (\ref{b4le}), (\ref{i2i3}) in (\ref{meb4}), we obtain that
\ce
\mE\(\sup_{0\leq t\leq T}|\hat{X}_{t}^{\e}-\bar{X}_{t}|^{2}\)
\leq C(\frac{\e}{\d}+\d)+C({\d}^{2}+{\d})+C\int_{0}^{T}\mE\(\sup_{0\leq r\leq s}|\hat{X}_{r}^{\e}-\bar{X}_{r}|^{2}\)\dif s.
\de
By the Gronwall inequality, it holds that
\ce
\mE\(\sup_{0\leq t\leq T}|\hat{X}_{t}^{\e}-\bar{X}_{t}|^{2}\)
\leq C\(\frac{\e}{\d}+{\d}^{2}+{\d}\).
\de

{\bf Step 2.} We prove (\ref{b4le}).

For any $t\in[0,T]$, it holds that
\be
I_1&=&6\Bigg(\mE\sup_{0\leq t\leq T}
 \Big|\int_{0}^{[\frac{t}{\d}]\d}\(b_{1}(X_{s(\d)}^{\e},\sL^{\mP}_{X_{s(\d)}^{\e}},\hat{Z}_{s}^{\e,z_0,\sL^{\mP}_{\xi}},\sL^{\mP}_{\hat{Z}_{s}^{\e,\xi}})
 -\bar{b}_{1}(X_{s(\d)}^{\e},\sL^{\mP}_{X_{s(\d)}^{\e}})\)\dif s\no\\
 &&\quad\quad\quad\quad+\int_{[\frac{t}{\d}]\d}^{t}\(b_{1}(X_{s(\d)}^{\e},\sL^{\mP}_{X_{s(\d)}^{\e}},\hat{Z}_{s}^{\e,z_0,\sL^{\mP}_{\xi}},\sL^{\mP}_{\hat{Z}_{s}^{\e,\xi}})
 -\bar{b}_{1}(X_{s(\d)}^{\e},\sL^{\mP}_{X_{s(\d)}^{\e}})\)\dif s\Big|^{2}\Bigg)\no\\
&\leq&12\Bigg(\mE\sup_{0\leq t\leq T}
\Big|\int_{0}^{[\frac{t}{\d}]\d}\(b_{1}(X_{s(\d)}^{\e},\sL^{\mP}_{X_{s(\d)}^{\e}},\hat{Z}_{s}^{\e,z_0,\sL^{\mP}_{\xi}},\sL^{\mP}_{\hat{Z}_{s}^{\e,\xi}})
-\bar{b}_{1}(X_{s(\d)}^{\e},\sL^{\mP}_{X_{s(\d)}^{\e}})\)\dif s\Big|^{2}\Bigg)\no\\
&&+12\Bigg(\mE\sup_{0\leq t\leq T}
\Big|\int_{[\frac{t}{\d}]\d}^{t}\(b_{1}(X_{s(\d)}^{\e},\sL^{\mP}_{X_{s(\d)}^{\e}},\hat{Z}_{s}^{\e,z_0,\sL^{\mP}_{\xi}},\sL^{\mP}_{\hat{Z}_{s}^{\e,\xi}})
-\bar{b}_{1}(X_{s(\d)}^{\e},\sL^{\mP}_{X_{s(\d)}^{\e}})\)\dif s\Big|^{2}\Bigg)\no\\
 &=:&B_{1}+B_{2}.
\label{b4cp}
\ee

Next, we estimate $B_{1}$. Note that
\be
B_{1}
&=&12\Bigg(\mE\sup_{0\leq t\leq T}
\Big|\int_{0}^{[\frac{t}{\d}]\d}\(b_{1}(X_{s(\d)}^{\e},\sL^{\mP}_{X_{s(\d)}^{\e}},\hat{Z}_{s}^{\e,z_0,\sL^{\mP}_{\xi}},\sL^{\mP}_{\hat{Z}_{s}^{\e,\xi}})
-\bar{b}_{1}(X_{s(\d)}^{\e},\sL^{\mP}_{X_{s(\d)}^{\e}})\)\dif s\Big|^{2}\Bigg)\no\\
&=&12\mE\Bigg(\sup_{0\leq t\leq T}
\Big|\sum\limits_{k=0}^{[\frac{t}{\d}]-1}\int_{k\d}^{(k+1)\d}\(b_{1}(X_{k\d}^{\e},\sL^{\mP}_{X_{k\d}^{\e}},\hat{Z}_{s}^{\e,z_0,\sL^{\mP}_{\xi}},\sL^{\mP}_{\hat{Z}_{s}^{\e,\xi}})
-\bar{b}_{1}(X_{k\d}^{\e},\sL^{\mP}_{X_{k\d}^{\e}})\)\dif s\Big|^{2}\Bigg)\no\\
&\leq&12\mE\Bigg(\sup_{0\leq t\leq T}[\frac{t}{\d}]\sum\limits_{k=0}^{[\frac{t}{\d}]-1}\left|\int_{k\d}^{(k+1)\d}\(b_{1}(X_{k\d}^{\e},\sL^{\mP}_{X_{k\d}^{\e}},\hat{Z}_{s}^{\e,z_0,\sL^{\mP}_{\xi}},\sL^{\mP}_{\hat{Z}_{s}^{\e,\xi}})
-\bar{b}_{1}(X_{k\d}^{\e},\sL^{\mP}_{X_{k\d}^{\e}})\)\dif s\right|^2\no\\
&\leq&12[\frac{T}{\d}]\sum\limits_{k=0}^{[\frac{T}{\d}]-1}
\mE\Bigg(\Big|\int_{k\d}^{(k+1)\d}\(b_{1}(X_{k\d}^{\e},\sL^{\mP}_{X_{k\d}^{\e}},\hat{Z}_{s}^{\e,z_0,\sL^{\mP}_{\xi}},\sL^{\mP}_{\hat{Z}_{s}^{\e,\xi}})
-\bar{b}_{1}(X_{k\d}^{\e},\sL^{\mP}_{X_{k\d}^{\e}})\)\dif s\Big|^{2}\Bigg)\no\\
&\leq&12[\frac{T}{\d}]^{2}\sup_{0\leq k\leq [\frac{T}{\d}]-1}
\mE\Bigg(\Big|\int_{k\d}^{(k+1)\d}\(b_{1}(X_{k\d}^{\e},\sL^{\mP}_{X_{k\d}^{\e}},\hat{Z}_{s}^{\e,z_0,\sL^{\mP}_{\xi}},\sL^{\mP}_{\hat{Z}_{s}^{\e,\xi}})
-\bar{b}_{1}(X_{k\d}^{\e},\sL^{\mP}_{X_{k\d}^{\e}})\)\dif s\Big|^{2}\Bigg)\no\\
&\leq&12\e^2(\frac{T}{\d})^{2}\sup_{0\leq k\leq [\frac{T}{\d}]-1}
\mE\Bigg(\Big|\int_{0}^{\d/\e}\(b_{1}(X_{k\d}^{\e},\sL^{\mP}_{X_{k\d}^{\e}},\hat{Z}_{\e s+k\d}^{\e,z_0,\sL^{\mP}_{\xi}},\sL^{\mP}_{\hat{Z}_{\e s+k\d}^{\e,\xi}})
-\bar{b}_{1}(X_{k\d}^{\e},\sL^{\mP}_{X_{k\d}^{\e}})\)\dif s\Big|^{2}\Bigg)\no\\
&=&12\e^2(\frac{T}{\d})^{2}\sup_{0\leq k\leq [\frac{T}{\d}]-1}\mE\Bigg<\int_{0}^{\d/\e}\(b_{1}(X_{k\d}^{\e},\sL^{\mP}_{X_{k\d}^{\e}},\hat{Z}_{\e s+k\d}^{\e,z_0,\sL^{\mP}_{\xi}},\sL^{\mP}_{\hat{Z}_{\e s+k\d}^{\e,\xi}})
-\bar{b}_{1}(X_{k\d}^{\e},\sL^{\mP}_{X_{k\d}^{\e}})\)\dif s,\no\\
&&\qquad\qquad \int_{0}^{\d/\e}\(b_{1}(X_{k\d}^{\e},\sL^{\mP}_{X_{k\d}^{\e}},\hat{Z}_{\e r+k\d}^{\e,z_0,\sL^{\mP}_{\xi}},\sL^{\mP}_{\hat{Z}_{\e r+k\d}^{\e,\xi}})
-\bar{b}_{1}(X_{k\d}^{\e},\sL^{\mP}_{X_{k\d}^{\e}})\)\dif r\Bigg>\no\\
&=&24\e^2(\frac{T}{\d})^{2}\sup_{0\leq k\leq [\frac{T}{\d}]-1}\int_{0}^{\d/\e}\int_r^{\d/\e}\mE\Bigg<\(b_{1}(X_{k\d}^{\e},\sL^{\mP}_{X_{k\d}^{\e}},\hat{Z}_{\e s+k\d}^{\e,z_0,\sL^{\mP}_{\xi}},\sL^{\mP}_{\hat{Z}_{\e s+k\d}^{\e,\xi}})
-\bar{b}_{1}(X_{k\d}^{\e},\sL^{\mP}_{X_{k\d}^{\e}})\),\no\\
&&\qquad\qquad \(b_{1}(X_{k\d}^{\e},\sL^{\mP}_{X_{k\d}^{\e}},\hat{Z}_{\e r+k\d}^{\e,z_0,\sL^{\mP}_{\xi}},\sL^{\mP}_{\hat{Z}_{\e r+k\d}^{\e,\xi}})
-\bar{b}_{1}(X_{k\d}^{\e},\sL^{\mP}_{X_{k\d}^{\e}})\)\Bigg>\dif s\dif r.
\label{b41c}
\ee

In the following, for $0<r<s\leq\d/\e$, set
\ce
\Phi(s,r)&:=&\mE\Bigg<b_{1}(X_{k\d}^{\e},\sL^{\mP}_{X_{k\d}^{\e}},\hat{Z}_{\e s+k\d}^{\e,z_0,\sL^{\mP}_{\xi}},\sL^{\mP}_{\hat{Z}_{\e s+k\d}^{\e,\xi}})-\bar{b}_{1}(X_{k\d}^{\e},\sL^{\mP}_{X_{k\d}^{\e}}),\\
&&\quad b_{1}(X_{k\d}^{\e},\sL^{\mP}_{X_{k\d}^{\e}},\hat{Z}_{\e r+k\d}^{\e,z_0,\sL^{\mP}_{\xi}},\sL^{\mP}_{\hat{Z}_{\e r+k\d}^{\e,\xi}})-\bar{b}_{1}(X_{k\d}^{\e},\sL^{\mP}_{X_{k\d}^{\e}})\Bigg>,
\de
and we estimate $\Phi(s,r)$. For any $s>0, \vartheta\in L^2(\Omega,\mathscr{F}_s,\mP,\mR^m), \varsigma\in L^2(\Omega,\mathscr{F}_s,\mP,\mR^n), \mu\in\cP_2(\mR^n), z\in\mR^m$, we consider two following equations
\ce
\check{Z}_t^{\e, s, \varsigma,\mu,\vartheta}&=&\vartheta+\frac{1}{\e} \int_s^t b_2(\varsigma, \mu, \check{Z}_r^{\e, s, \varsigma,\mu,\vartheta},\sL^{\mP}_{\check{Z}_r^{\e, s, \varsigma,\mu,\vartheta}})\dif r\\
&&+\frac{1}{\sqrt{\e}} \int_s^t \s_2(\varsigma, \mu, \check{Z}_r^{\e, s, \varsigma,\mu,\vartheta},\sL^{\mP}_{\check{Z}_r^{\e, s, \varsigma,\mu,\vartheta}})\dif W_r, \\
\check{Z}_t^{\e, s, \varsigma,\mu,z,\sL^\mP_\vartheta}&=&z+\frac{1}{\e} \int_s^t b_2(\varsigma, \mu, \check{Z}_r^{\e, s, \varsigma,\mu,z,\sL^\mP_\vartheta},\sL^{\mP}_{\check{Z}_r^{\e, s, \varsigma,\mu,\vartheta}})\dif r\\
&&+\frac{1}{\sqrt{\e}} \int_s^t \s_2(\varsigma, \mu, \check{Z}_r^{\e, s, \varsigma,\mu,z,\sL^\mP_\vartheta},\sL^{\mP}_{\check{Z}_r^{\e, s, \varsigma,\mu,\vartheta}})\dif W_r,\quad t \geqslant s.
\de
Then it holds that
$$
\hat{Z}_t^{\e,\xi}=\check{Z}_t^{\e, k \delta, X_{k \delta}^\e, \mathscr{L}^\mP_{X_{k \delta}^\e},\hat{Z}_{k \delta}^{\e,\xi}}, \quad \hat{Z}_t^{\e,z_0,\sL_\xi}=\check{Z}_t^{\e, k \delta, X_{k \delta}^\e, \mathscr{L}^\mP_{X_{k \delta}^\e}, \hat{Z}_{k \delta}^{\e,z_0,\sL^\mP_\xi},\sL^\mP_{\hat{Z}_{k \delta}^{\e,\xi}}}\quad t \in[k \delta,(k+1) \delta],
$$
and furthermore
\ce
\Phi(s,r)&=&\mE\Bigg<b_{1}(X_{k\d}^{\e},\sL^{\mP}_{X_{k\d}^{\e}},\check{Z}^{\e, k \delta, X_{k \delta}^\e, \mathscr{L}^\mP_{X_{k \delta}^\e}, \hat{Z}_{k \delta}^{\e,z_0,\sL^\mP_\xi},\sL^\mP_{\hat{Z}_{k \delta}^{\e,\xi}}}_{\e s+k\d},\sL^\mP_{\check{Z}_{\e s+k\d}^{\e, k \delta, X_{k \delta}^\e, \mathscr{L}^\mP_{X_{k \delta}^\e}, \hat{Z}_{k \delta}^{\e,\xi}}})
-\bar{b}_{1}(X_{k\d}^{\e},\sL^{\mP}_{X_{k\d}^{\e}}),\\
&&b_{1}(X_{k\d}^{\e},\sL^{\mP}_{X_{k\d}^{\e}},\check{Z}^{\e, k \delta, X_{k \delta}^\e, \mathscr{L}^\mP_{X_{k \delta}^\e}, \hat{Z}_{k \delta}^{\e,z_0,\sL^\mP_\xi},\sL^\mP_{\hat{Z}_{k \delta}^{\e,\xi}}}_{\e r+k\d},\sL^\mP_{\check{Z}_{\e r+k\d}^{\e, k \delta, X_{k \delta}^\e, \mathscr{L}^\mP_{X_{k \delta}^\e}, \hat{Z}_{k \delta}^{\e,\xi}}})
-\bar{b}_{1}(X_{k\d}^{\e},\sL^{\mP}_{X_{k\d}^{\e}})\Bigg>.
\de
Note that $X_{k \delta}^\e, \hat{Z}_{k \delta}^{\e,z_0,\sL^\mP_\xi}$ are $\sF_{k\d}$-measurable, and for any $x\in\mR^n$, $\check{Z}_t^{\e, k\d, x,\sL^{\mP}_{X_{k\d}^{\e}},z,\sL^\mP_{\hat{Z}_{k \delta}^{\e,\xi}}}$ is independent of $\sF_{k\d}$. Thus, we have that
\ce
\Phi(s,r)&=&\mE\Bigg[\mE\Bigg[\Bigg<b_{1}(X_{k\d}^{\e},\sL^{\mP}_{X_{k\d}^{\e}},\check{Z}^{\e, k \delta, X_{k \delta}^\e, \mathscr{L}^\mP_{X_{k \delta}^\e}, \hat{Z}_{k \delta}^{\e,z_0,\sL^\mP_\xi},\sL^\mP_{\hat{Z}_{k \delta}^{\e,\xi}}}_{\e s+k\d},\sL^\mP_{\check{Z}_{\e s+k\d}^{\e, k \delta, X_{k \delta}^\e, \mathscr{L}^\mP_{X_{k \delta}^\e}, \hat{Z}_{k \delta}^{\e,\xi}}})
-\bar{b}_{1}(X_{k\d}^{\e},\sL^{\mP}_{X_{k\d}^{\e}}),\\
&&b_{1}(X_{k\d}^{\e},\sL^{\mP}_{X_{k\d}^{\e}},\check{Z}^{\e, k \delta, X_{k \delta}^\e, \mathscr{L}^\mP_{X_{k \delta}^\e}, \hat{Z}_{k \delta}^{\e,z_0,\sL^\mP_\xi},\sL^\mP_{\hat{Z}_{k \delta}^{\e,\xi}}}_{\e r+k\d},\sL^\mP_{\check{Z}_{\e r+k\d}^{\e, k \delta, X_{k \delta}^\e, \mathscr{L}^\mP_{X_{k \delta}^\e}, \hat{Z}_{k \delta}^{\e,\xi}}})
-\bar{b}_{1}(X_{k\d}^{\e},\sL^{\mP}_{X_{k\d}^{\e}})\Bigg>\Bigg{|}\sF_{k\d}\Bigg]\Bigg]\\
&=&\mE\Bigg[\Bigg(\mE\Bigg<b_{1}(x,\mu,\check{Z}^{\e, k \delta, x, \mu, z,\nu}_{\e s+k\d},\sL^\mP_{\check{Z}_{\e r+k\d}^{\e, k \delta, x, \mu, \hat{Z}_{k \delta}^{\e,\xi}}})-\bar{b}_{1}(x,\mu),\\
&&\qquad\qquad b_{1}(x,\mu,\check{Z}^{\e, k \delta, x, \mu, z,\nu}_{\e r+k\d},\sL^\mP_{\check{Z}_{\e r+k\d}^{\e, k \delta, x, \mu, \hat{Z}_{k \delta}^{\e,\xi}}})
-\bar{b}_{1}(x,\mu)\Bigg>\Bigg)\Bigg{|}_{(x,\mu,z,\nu)=(X_{k \delta}^\e,\sL^{\mP}_{X_{k\d}^{\e}},\hat{Z}_{k \delta}^{\e,z_0,\sL_\xi},\sL_{\hat{Z}_{k \delta}^{\e,\xi}})}\Bigg].
\de

Here, we investigate $\check{Z}^{\e, k \delta, x, \mu, z,\nu}_{\e s+k\d}$. On one hand, it holds that
\ce
\check{Z}^{\e, k \delta, x, \mu, z,\nu}_{\e s+k\d}&=&z+\frac{1}{\e} \int_{k\d}^{\e s+k\d}b_2(x, \mu, \check{Z}^{\e, k \delta, x, \mu, z,\nu}_{r},\sL^{\mP}_{\check{Z}^{\e, k \delta, x, \mu, \hat{Z}_{k \delta}^{\e,\xi}}_{r}})\dif r\\
&&+\frac{1}{\sqrt{\e}} \int_{k\d}^{\e s+k\d} \s_2(x, \mu, \check{Z}^{\e, k \delta, x, \mu, z,\nu}_{r},\sL^{\mP}_{\check{Z}^{\e, k \delta, x, \mu, \hat{Z}_{k \delta}^{\e,\xi}}_{r}})\dif W_r\\
&=&z+\frac{1}{\e} \int_{0}^{\e s}b_2(x, \mu, \check{Z}^{\e, k \delta, x, \mu, z,\nu}_{u+k\d},\sL^{\mP}_{\check{Z}^{\e, k \delta, x, \mu, \hat{Z}_{k \delta}^{\e,\xi}}_{u+k\d}})\dif u\\
&&+\frac{1}{\sqrt{\e}} \int_{0}^{\e s} \s_2(x, \mu, \check{Z}^{\e, k \delta, x, \mu, z,\nu}_{u+k\d},\sL^{\mP}_{\check{Z}^{\e, k \delta, x, \mu, \hat{Z}_{k \delta}^{\e,\xi}}_{u+k\d}})\dif \check{W}_u\\
&=&z+\int_{0}^{s}b_2(x, \mu, \check{Z}^{\e, k \delta, x, \mu, z,\nu}_{\e v+k\d},\sL^{\mP}_{\check{Z}^{\e, k \delta, x, \mu, \hat{Z}_{k \delta}^{\e,\xi}}_{\e v+k\d}})\dif v\\
&&+\int_{0}^{s} \s_2(x, \mu, \check{Z}^{\e, k \delta, x, \mu, z,\nu}_{\e v+k\d},\sL^{\mP}_{\check{Z}^{\e, k \delta, x, \mu,\hat{Z}_{k \delta}^{\e,\xi}}_{\e v+k\d}})\dif \check{\check{W}}_v,
\de
where $\check{W}_u:=W_{u+k\d}-W_{k\d}$ and $\check{\check{W}}_v:=\frac{1}{\sqrt{\e}}\check{W}_{\e v}$ are two $m$-dimensional standard Brownian motions. On the other hand, we notice that Eq.(\ref{Eq21}) is just written as
\ce
Z_{s}^{x,\mu,z,\nu}=z+\int_0^s b_{2}(x,\mu,Z_{r}^{x,\mu,z,\nu},\sL^{\mP}_{Z_{r}^{x,\mu,\zeta}})\dif r+\int_0^s\s_{2}(x,\mu,Z_{r}^{x,\mu,z,\nu},\sL^{\mP}_{Z_{r}^{x,\mu,\zeta}})\dif W_{r},
\de
where $\zeta\in L^{2}(\Omega,\sF_0,\mP;\mR^m), \sL^{\mP}_{\zeta}=\nu$. That is, for $s\in[0,\d/\e]$, $(\check{Z}^{\e, k \delta, x, \mu, z,\nu}_{\e s+k\d},\sL^{\mP}_{\check{Z}^{\e, k \delta, x, \mu, \hat{Z}_{k \delta}^{\e,\xi}}_{\e s+k\d}})$ and $(Z_{s}^{x,\mu,z,\nu},\sL^{\mP}_{Z_{s}^{x,\mu,\zeta}})$ have the same distribution. But $(Z_{\cdot}^{x,\mu,z,\nu},\sL^{\mP}_{Z_{\cdot}^{x,\mu,\zeta}})$ is {\it not} a Markov process. Therefore, we need to construct a Markov process based on $(Z_{\cdot}^{x,\mu,z,\nu},\sL^{\mP}_{Z_{\cdot}^{x,\mu,\zeta}})$.

Let $C([0,\infty),\mR^m)$ be the collection of continuous functions from $[0,\infty)$ to $\mR^m$ with the uniform convergence topology. Set
\ce
&&\tilde{\Omega}:=C([0,\infty),\mR^m)\times C([0,\infty),\cP_2(\mR^m)), \\
&&\tilde{\sF}=\sB\(C([0,\infty),\mR^m)\)\times \sB\(C\left([0,\infty),\cP_2(\mR^m)\right)\),\\
&&\tilde{\sF}_t=\sigma(M_{r}, 0\leq r\leq t), \quad t\geq 0,
\de
where $M_{\cdot}$ is the coordinate process. Then by \cite[Theorem 4.11]{rrw}, there exists a unique probability measure $\tilde{\mP}^{x,\mu,z,\nu}$ on $(\tilde{\Omega}, \tilde{\sF})$ such that $M_{\cdot}$ is a Markov process with respect to $(\tilde{\sF}_t)$ with the transition function $\{{\bf P}^{x,\mu}_{t}(z,\nu;\cdot)=\sL_{Z^{x,\mu, z,\nu}_{t}}\times\delta_{\sL^{\mP}_{Z_{t}^{x, \mu, \zeta}}}: t\geq 0, (z,\nu)\in\mR^m\times\cP_2(\mR^m)\}$ and $\sL^{\tilde{\mP}^{x,\mu,z,\nu}}_{M_{0}}=\d_z\times\d_\nu$. Note that
\ce
\sL^{\tilde{\mP}^{x,\mu,z,\nu}}_{M_{t}}:=\tilde{\mP}^{x,\mu,z,\nu}\circ M^{-1}_{t}=\int_{\mR^m\times\cP_2(\mR^m)}{\bf P}^{x,\mu}_{t}(z',\nu';\cdot)\d_z\times\d_\nu(\dif z', \dif \nu')=\sL_{Z^{x,\mu,z,\nu}_{t}}\times\delta_{\sL^{\mP}_{Z_{t}^{x, \mu,\zeta}}}.
\de
Thus, set 
$$
\tilde{\mP}:=\tilde{\mP}^{x,\mu,z,\nu}, \quad (\tilde{Z}^{x,\mu,z,\nu}_{t},\sL^{\tilde{\mP}}_{\tilde{Z}_{t}^{x, \mu, \tilde{\zeta}}}):=M_{t}, 
$$
and $\sL^{\tilde{\mP}}_{\tilde{Z}^{x,\mu,z,\nu}_{t}}=\sL^{\mP}_{Z^{x,\mu,z,\nu}_{t}}, \sL^{\tilde{\mP}}_{\tilde{Z}_{t}^{x, \mu, \tilde{\zeta}}}=\sL^{\mP}_{Z_{t}^{x, \mu,\zeta}},\sL^{\tilde{\mP}}_{\tilde{\zeta}}=\sL^{\mP}_{\zeta}=\nu$, which implies that
\ce
&&\mE\Bigg<b_{1}(x,\mu,\check{Z}^{\e, k \delta, x, \mu, z,\nu}_{\e s+k\d},\sL_{\check{Z}_{\e s+k\d}^{\e, k \delta, x, \mu, \hat{Z}_{k \delta}^{\e,\xi}}})-\bar{b}_{1}(x,\mu),b_{1}(x,\mu,\check{Z}^{\e, k \delta, x, \mu, z,\nu}_{\e r+k\d},\sL_{\check{Z}_{\e r+k\d}^{\e, k \delta, x, \mu, \hat{Z}_{k \delta}^{\e,\xi}}})
-\bar{b}_{1}(x,\mu)\Bigg>\\
&=&\tilde{\mE}\Bigg<b_{1}(x,\mu,\tilde{Z}_{s}^{x,\mu,z,\nu},\sL^{\tilde{\mP}}_{\tilde{Z}_{s}^{x,\mu,\tilde{\zeta}}})-\bar{b}_{1}(x,\mu),b_{1}(x,\mu,\tilde{Z}_{r}^{x,\mu,z,\nu},\sL^{\tilde{\mP}}_{\tilde{Z}_{r}^{x,\mu,\tilde{\zeta}}})
-\bar{b}_{1}(x,\mu)\Bigg>\\
&=&\tilde{\mE}\left[\tilde{\mE}\left[\Bigg<b_{1}(x,\mu,\tilde{Z}_{s}^{x,\mu,z,\nu},\sL^{\tilde{\mP}}_{\tilde{Z}_{s}^{x,\mu,\tilde{\zeta}}})-\bar{b}_{1}(x,\mu), b_{1}(x,\mu,\tilde{Z}_{r}^{x,\mu,z,\nu},\sL^{\tilde{\mP}}_{\tilde{Z}_{r}^{x,\mu,\tilde{\zeta}}})
-\bar{b}_{1}(x,\mu)\Bigg>\Bigg{|}\tilde{\sF}_r\right]\right]\\
&=&\tilde{\mE}\left[\Bigg<\tilde{\mE}\left[b_{1}(x,\mu,\tilde{Z}_{s}^{x,\mu,z,\nu},\sL^{\tilde{\mP}}_{\tilde{Z}_{s}^{x,\mu,\tilde{\zeta}}})\Bigg{|}\tilde{\sF}_r\right]-\bar{b}_{1}(x,\mu),b_{1}(x,\mu,\tilde{Z}_{r}^{x,\mu,z,\nu},\sL^{\tilde{\mP}}_{\tilde{Z}_{r}^{x,\mu,\tilde{\zeta}}})
-\bar{b}_{1}(x,\mu)\Bigg>\right]\\
&\leq&\left(\tilde{\mE}\left|\tilde{\mE}\left[b_{1}(x,\mu,\tilde{Z}_{s}^{x,\mu,z,\nu},\sL^{\tilde{\mP}}_{\tilde{Z}_{s}^{x,\mu,\tilde{\zeta}}})\Bigg{|}\tilde{\sF}_r\right]-\bar{b}_{1}(x,\mu)\right|^2\right)^{1/2}\left(\tilde{\mE}|b_{1}(x,\mu,\tilde{Z}_{r}^{x,\mu,z,\nu},\sL^{\tilde{\mP}}_{\tilde{Z}_{r}^{x,\mu,\tilde{\zeta}}})
-\bar{b}_{1}(x,\mu)|^2\right)^{1/2}.
\de
Moreover, based on (\ref{meu2}), (\ref{memu2}) and (\ref{zxmuxi}), we obtain that
\ce
&&\left(\tilde{\mE}\left|\tilde{\mE}\left[b_{1}(x,\mu,\tilde{Z}_{s}^{x,\mu,z,\nu},\sL^{\tilde{\mP}}_{\tilde{Z}_{s}^{x,\mu,\tilde{\zeta}}})\Bigg{|}\tilde{\sF}_r\right]-\bar{b}_{1}(x,\mu)\right|^2\right)^{1/2}\\
&=&\left(\tilde{\mE}\left|\tilde{\mE}\left[b_{1}(x,\mu,\tilde{Z}_{s-r}^{x,\mu,\hat{z},\hat{\nu}},\sL^{\tilde{\mP}}_{\tilde{Z}_{s-r}^{x,\mu,\tilde{Z}_{r}^{x,\mu,\tilde{\zeta}}}})\right]\Bigg{|}_{(\hat{z},\hat{\nu})=(\tilde{Z}_{r}^{x,\mu,z,\nu},\sL^{\tilde{\mP}}_{\tilde{Z}_{r}^{x,\mu,\tilde{\zeta}}})}-\bar{b}_{1}(x,\mu)\right|^2\right)^{1/2}\\
&\leq&Ce^{-(\b_{1}-L_{b_2,\s_2})(s-r)/2}\left(\tilde{\mE}|\tilde{Z}_{r}^{x,\mu,\tilde{\zeta}}|^2+1+|x|^2+\|\mu\|^2+\tilde{\mE}|\tilde{Z}_{r}^{x,\mu,z,\nu}|^2\right)^{1/2}\\
&\leq&Ce^{-(\b_{1}-L_{b_2,\s_2})(s-r)/2}(1+|x|+\|\mu\|+|z|+\|\nu\|).
\de
and
\ce
&&\left(\tilde{\mE}|b_{1}(x,\mu,\tilde{Z}_{r}^{x,\mu,z,\nu},\sL^{\tilde{\mP}}_{\tilde{Z}_{r}^{x,\mu,\tilde{\zeta}}})-\bar{b}_{1}(x,\mu)|^2\right)^{1/2}\\
&\leq&\(C(\tilde{\mE}|\tilde{Z}_{r}^{x,\mu,\tilde{\zeta}}|^2+1+|x|^2+\|\mu\|^2+\tilde{\mE}|\tilde{Z}_{r}^{x,\mu,z,\nu}|^2)\)^{\frac{1}{2}}\\
&\leq&C(1+|x|+\|\mu\|+|z|+\|\nu\|).
\de

Combining the above deduction, by Lemma \ref{xtztc}, \ref{hatze} one can have that
\ce
\Phi(s,r)\leq Ce^{-(\b_{1}-L_{b_2,\s_2})(s-r)/2}.
\de
Inserting the above inequality in (\ref{b41c}), we get that
\be
B_1\leq C(\frac{\e}{\d})^2\int_{0}^{\frac{\d}{\e}}\int_{r}^{\frac{\d}{\e}}Ce^{-(\b_{1}-L_{b_2,\s_2})(s-r)/2}\dif s\dif r\leq C\frac{\e}{\d}.
\label{b4de}
\ee

Next, we estimate $B_{2}$. By (\ref{b1line}), Lemma \ref{xtztc}, \ref{hatze} and the H\"older inequality, one could get that
\ce
B_{2}&\leq& 24\d\mE\sup_{0\leq t\leq T}
\int_{[\frac{t}{\d}]\d}^{t}\(|b_{1}(X_{s(\d)}^{\e},\sL^{\mP}_{X_{s(\d)}^{\e}},\hat{Z}_{s}^{\e,z_0,\sL^{\mP}_{\xi}},\sL^{\mP}_{\hat{Z}_{s}^{\e,\xi}})|^2
+|\bar{b}_{1}(X_{s(\d)}^{\e},\sL^{\mP}_{X_{s(\d)}^{\e}})|^2\)\dif s\\
&\leq&C\d\mE\int_0^T(1+|X_{s(\d)}^{\e}|^2+\|\sL^{\mP}_{X_{s(\d)}^{\e}}\|^2+|\hat{Z}_{s}^{\e,z_0,\sL^{\mP}_{\xi}}|^2+\|\sL^{\mP}_{\hat{Z}_{s}^{\e,\xi}}\|^2)\dif s\\
&\leq& C\d,
\de
which together with (\ref{b4de}) implies (\ref{b4le}). The proof is complete.
\end{proof}

At present, we are ready to prove Theorem \ref{xbarxp}.

{\bf Proof of Theorem \ref{xbarxp}.} Taking $\d=\e^\g$, by (\ref{xehatxe}) and Lemma \ref{2orde}, we get that 
$$
\mE\(\sup_{0\leq t\leq T}|X_{t}^{\e}-\bar{X}_{t}|^{2}\)\leq C\(\e^{1-\g}+\e^{2\g}+\e^\g\).
$$
This is just (\ref{2or}). 

Next, by the Chebyshev inequality and (\ref{2or}), it holds that for any $\t>0$
$$ 
\mP\(\sup_{0\leq t\leq T}|X_{t}^{\e}-\bar{X}_{t}|>\t\)\leq \frac{\mE\(\sup\limits_{0\leq t\leq T}|X_{t}^{\e}-\bar{X}_{t}|^{2}\)}{\t^2}\leq \frac{C}{\t^2}\(\e^{1-\g}+\e^{2\g}+\e^\g\),
$$
which implies that
$$
\sup_{0\leq t\leq T}|X_{t}^{\e}-\bar{X}_{t}|\overset{\mP}{\rightarrow} 0,
$$
as $\e$ tends to $0$. Besides, from Lemma \ref{xtztc} and (\ref{barxb}) it follows that 
$$
\sup\limits_{\e}\mE\sup\limits_{0\leq t\leq T}|X_{t}^{\e}-\bar{X}_{t}|^{2p+2}\leq C(1+\mE|\varrho|^{2p+2}+|z_0|^{2p+2}+\mE|\xi|^{2p+2}).
$$
Therefore, by the Vitali convergence theorem one can obtain that
$$
\lim\limits_{\e\rightarrow 0}\mE\(\sup_{0\leq t\leq T}|X_{t}^{\e}-\bar{X}_{t}|^{2p}\)=0,
$$
which completes the proof.

\section{Proof of Theorem \ref{effifilt}}\label{prooseco}

In this section, we prove Theorem \ref{effifilt}. First of all, we prepare an important lemma.

\bl\label{4}
Under the assumption $(\mathbf{H}_{h})$, there exists a constant $C>0$ such that
\ce
\mE|\rho_{t}^{0}(1)|^{-r}\leq \exp\{(2r^{2}+r+1)CT/2\}, \quad r>1.
\de
\el

Since its proof is similar to that of \cite[Lemma 3.6]{Qiao1}, we omit it.

Now, we are ready to prove Theorem \ref{effifilt}.

{\bf Proof of Theorem \ref{effifilt}.} 
{\bf Step 1.} We estimate $\pi_{t}^{\e}(F)-\pi_{t}^{0}(F)$.

Based on the H\"{o}lder inequality and these definitions of $\pi_{t}^{\e}(F)$ and $\pi_{t}^{0}(F)$, we get that
\be
\mE|\pi_{t}^{\e}(F)-\pi_{t}^{0}(F)|^q
&=&\mE\left|\frac{\rho_{t}^{\e}(F)-\rho_{t}^{0}(F)}{\rho_{t}^{0}(1)}
-\pi_{t}^{\e}(F)\frac{\rho_{t}^{\e}(1)-\rho_{t}^{0}(1)}{\rho_{t}^{0}(1)}\right|^q\no\\
&\leq& 2^{q-1}\mE\left|\frac{\rho_{t}^{\e}(F)-\rho_{t}^{0}(F)}{\rho_{t}^{0}(1)}\right|^q
+2^{q-1}\mE\left|\pi_{t}^{\e}(F)\frac{\rho_{t}^{\e}(1)-\rho_{t}^{0}(1)}{\rho_{t}^{0}(1)}\right|^q\no\\
&\leq&2^{q-1}\left(\mE|\rho_{t}^{\e}(F)-\rho_{t}^{0}(F)|^{2q}\right)^{\frac{1}{2}}
\left(\mE|\rho_{t}^{0}(1)|^{-2q}\right)^{\frac{1}{2}}\no\\
&&+2^{q-1}\|F\|^q_{C_{b,lip}(\mR^n\times\cP_2(\mR^n))}\left(\mE|\rho_{t}^{\e}(1)-\rho_{t}^{0}(1)|^{2q}\right)^{\frac{1}{2}}
\left(\mE|\rho_{t}^{0}(1)|^{-2q}\right)^{\frac{1}{2}}\no\\
&\leq&C\left(\mE|\rho_{t}^{\e}(F)-\rho_{t}^{0}(F)|^{2q}\right)^{\frac{1}{2}}\no\\
&&+C\|F\|^q_{C_{b,lip}(\mR^n\times\cP_2(\mR^n))}\left(\mE|\rho_{t}^{\e}(1)-\rho_{t}^{0}(1)|^{2q}\right)^{\frac{1}{2}},
\label{9}
\ee
where $\|F\|_{C_{b,lip}(\mR^n\times\cP_2(\mR^n))}$ denotes the norm of $F$ in $C_{b,lip}(\mR^n\times\cP_2(\mR^n))$.

By the deduction in {\bf Step 2}, it holds that
\be
\lim_{\e\rightarrow0}\mE^{\mP^{\e}}[|\rho_{t}^{\e}(F)-\rho_{t}^{0}(F)|^{2q}]=0.
\label{rho}
\ee
Inserting (\ref{rho}) in (\ref{9}), we obtain that
$$
\lim_{\e\rightarrow0}\mE|\pi_{t}^{\e}(F)-\pi_{t}^{0}(F)|^q=0.
$$

{\bf Step 2.} We prove (\ref{rho}).

By the measure transformation and the H\"{o}lder inequality, it holds that
\be
\mE|\rho_{t}^{\e}(F)-\rho_{t}^{0}(F)|^{2q}
=\mE^{\mP^{\e}}[|\rho_{t}^{\e}(F)-\rho_{t}^{0}(F)|^{2q}\Lambda_{T}^{\e}]
\leq \left(\mE^{\mP^{\e}}|\rho_{t}^{\e}(F)-\rho_{t}^{0}(F)|^{4q}\right)^{\frac{1}{2}}
\(\mE^{\mP^{\e}}(\Lambda_{T}^{\e})^{2}\)^{\frac{1}{2}}.
\label{10}
\ee
On one hand, it is not difficult to prove that
$$
(\mE^{\mP^{\e}}(\Lambda_{T}^{\e})^{2})^{\frac{1}{2}}\leq \exp\{CT\}.
$$
On the other hand, from the Jensen inequality, it follows that
\be
&&\mE^{\mP^{\e}}[|\rho_{t}^{\e}(F)-\rho_{t}^{0}(F)|^{4q}]\no\\
&=&\mE^{\mP^{\e}}\left[\left|\mE^{\mP^{\e}}\left(F(X_{t}^{\e},\sL^{\mP}_{X_{t}^{\e}})\Lambda_{t}^{\e}|\mathscr{F}_{t}^{Y^{\e}}\right)-
\mE^{\mP^{\e}}\left(F(\bar{X}_{t},\sL^{\mP}_{\bar{X}_{t}})\Lambda_{t}^{0}|\mathscr{F}_{t}^{Y^{\e}}\right)\right|^{4q}\right]\no\\
&=&\mE^{\mP^{\e}}\left[\left|\mE^{\mP^{\e}}\left(F(X_{t}^{\e},\sL^{\mP}_{X_{t}^{\e}})\Lambda_{t}^{\e}
-F(\bar{X}_{t},\sL^{\mP}_{\bar{X}_{t}})\Lambda_{t}^{0}|\mathscr{F}_{t}^{Y^{\e}}\right)\right|^{4q}\right]\no\\
&\leq&\mE^{\mP^{\e}}\left[\mE^{\mP^{\e}}\left(|F(X_{t}^{\e},\sL^{\mP}_{X_{t}^{\e}})\Lambda_{t}^{\e}
-F(\bar{X}_{t},\sL^{\mP}_{\bar{X}_{t}})\Lambda_{t}^{0}|^{4q}\Big|\mathscr{F}_{t}^{Y^{\e}}\right)\right]\no\\
&=&\mE^{\mP^{\e}}|F(X_{t}^{\e},\sL^{\mP}_{X_{t}^{\e}})\Lambda_{t}^{\e}
-F(\bar{X}_{t},\sL^{\mP}_{\bar{X}_{t}})\Lambda_{t}^{0}|^{4q}\no\\
&\leq&2^{4q-1}\mE^{\mP^{\e}}|F(X_{t}^{\e},\sL^{\mP}_{X_{t}^{\e}})\Lambda_{t}^{\e}
-F(\bar{X}_{t},\sL^{\mP}_{\bar{X}_{t}})\Lambda_{t}^{\e}|^{4q}\no\\
&&+2^{4q-1}\mE^{\mP^{\e}}|F(\bar{X}_{t},\sL^{\mP}_{\bar{X}_{t}})\Lambda_{t}^{\e}
-F(\bar{X}_{t},\sL^{\mP}_{\bar{X}_{t}})\Lambda_{t}^{0}|^{4q}\no\\
&=:&A_{1}+A_{2}.
\label{11}
\ee

For $A_{1}$, by the H\"{o}lder inequality, we know that
\ce
A_{1}
&=&2^{4q-1}\mE^{\mP^{\e}}|F(X_{t}^{\e},\sL^{\mP}_{X_{t}^{\e}})\Lambda_{t}^{\e}
-F(\bar{X}_{t},\sL^{\mP}_{\bar{X}_{t}})\Lambda_{t}^{\e}|^{4q}\\
&\leq&2^{4q-1}\(\mE^{\mP^{\e}}|F(X_{t}^{\e},\sL^{\mP}_{X_{t}^{\e}})
-F(\bar{X}_{t},\sL^{\mP}_{\bar{X}_{t}})|^{8q}\)^{\frac{1}{2}}
\(\mE^{\mP^{\e}}|\Lambda_{t}^{\e}|^{8q}\)^{\frac{1}{2}}\\
&\leq&2^{4q-1}\Big(\|F\|_{C_{b,lip}(\mR^n\times\cP_2(\mR^n))}^{8q}\mE^{\mP^{\e}}(|X_{t}^{\e}-\bar{X}_{t}|+\mW_2(\sL^{\mP}_{X_{t}^{\e}},\sL^{\mP}_{\bar{X}_{t}}))^{8q}\)^{\frac{1}{2}}
\(\mE^{\mP^{\e}}|\Lambda_{t}^{\e}|^{8q}\)^{\frac{1}{2}}\\
&\leq& 2^{4q-1}\|F\|_{C_{b,lip}(\mR^n\times\cP_2(\mR^n))}^{4q}
\(\mE^{\mP^{\e}}\(|X_{t}^{\e}-\bar{X}_{t}|+(\mE|X_{t}^{\e}-\bar{X}_{t}|^2)^{1/2}\)^{8q}\)^{\frac{1}{2}}
\(\mE^{\mP^{\e}}|\Lambda_{t}^{\e}|^{8q}\)^{\frac{1}{2}}\\
&\leq& 2^{4q-1}\|F\|_{C_{b,lip}(\mR^n\times\cP_2(\mR^n))}^{4q}
\(\mE^{\mP^{\e}}|\Lambda_{t}^{\e}|^{8q}\)^{\frac{1}{2}}\\
&&\cdot\(\mE^{\mP^{\e}}\(2^{8q-1}|X_{t}^{\e}-\bar{X}_{t}|^{8q}
+2^{8q-1}(\mE|X_{t}^{\e}-\bar{X}_{t}|^2)^{4q}\)\)^{\frac{1}{2}}\\
&\leq&2^{4q-1}\|F\|_{C_{b,lip}(\mR^n\times\cP_2(\mR^n))}^{4q}
\(\mE^{\mP^{\e}}|\Lambda_{t}^{\e}|^{8q}\)^{\frac{1}{2}}\\
&&\cdot\(\mE^{\mP^{\e}}\(2^{8q-1}|X_{t}^{\e}-\bar{X}_{t}|^{8q}+2^{8q-1}
\mE|X_{t}^{\e}-\bar{X}_{t}|^{8q}\)\)^{\frac{1}{2}}\\
&\leq&2^{8q-1-\frac{1}{2}}\|F\|_{C_{b,lip}(\mR^n\times\cP_2(\mR^n))}^{4q}
\(\mE^{\mP^{\e}}|\Lambda_{t}^{\e}|^{8q}\)^{\frac{1}{2}}\\
&&\cdot\(\mE^{\mP^{\e}}\(|X_{t}^{\e}-\bar{X}_{t}|^{8q}\)
+\mE|X_{t}^{\e}-\bar{X}_{t}|^{8q}\)^{\frac{1}{2}}.
\de
Then we estimate $\mE^{\mP^{\e}}(|X_{t}^{\e}-\bar{X}_{t}|^{8q})+\mE|X_{t}^{\e}-\bar{X}_{t}|^{8q}$. Note that
\ce
\mE^{\mP^{\e}}|X_{t}^{\e}-\bar{X}_{t}|^{8q}+\mE|X_{t}^{\e}-\bar{X}_{t}|^{8q}
&=&\mE[|X_{t}^{\e}-\bar{X}_{t}|^{8q}(\Lambda_{T}^{\e})^{-1}]+\mE|X_{t}^{\e}-\bar{X}_{t}|^{8q}\\
&\leq&
\(\mE|X_{t}^{\e}-\bar{X}_{t}|^{16q}\)^{\frac{{1}}{2}}
\(\mE(\Lambda_{T}^{\e})^{-2}\)^{\frac{1}{2}}+
\mE|X_{t}^{\e}-\bar{X}_{t}|^{8q}.
\de
 By the similar deduction to that in Lemma \ref{4}, one can get that
$$
\mE(\Lambda_{T}^{\e})^{-2}\leq C,
$$
which yields that
\be
A_{1}\leq 2^{8q-\frac{3}{2}}\|F\|_{C_{b,lip}(\mR^n\times\cP_2(\mR^n))}^{4q}\left[C\(\mE|X_{t}^{\e}-\bar{X}_{t}|^{16q}\)^{\frac{{1}}{2}}+
\mE|X_{t}^{\e}-\bar{X}_{t}|^{8q}\right].
\label{12}
\ee
Combining (\ref{mesu}) with (\ref{12}) and taking the limit as $\e\rightarrow0$, we get that
\be
\lim_{\e\rightarrow0}A_{1}=0.
\label{13}
\ee

For $A_{2}$, by the H\"{o}lder inequality, one can obtain that
\ce
A_{2}
&=&2^{4q-1}\mE^{\mP^{\e}}\left|F(\bar{X}_{t},\sL^{\mP}_{\bar{X}_{t}})\Lambda_{t}^{\e}
-F(\bar{X}_{t},\sL^{\mP}_{\bar{X}_{t}})\Lambda_{t}^{0}\right|^{4q}\\
&\leq&2^{4q-1}\|F\|_{C_{b,lip}(\mR^n\times\cP_2(\mR^n))}^{4q}
\mE^{\mP^{\e}}|\Lambda_{t}^{\e}-\Lambda_{t}^{0}|^{4q}\\
&=&2^{4q-1}\|F\|_{C_{b,lip}(\mR^n\times\cP_2(\mR^n))}^{4q}
\mE\Big[|\Lambda_{t}^{\e}-\Lambda_{t}^{0}|^{4q}(\Lambda_{T}^{\e})^{-1}\Big]\\
&\leq&2^{4q-1}\|F\|_{C_{b,lip}(\mR^n\times\cP_2(\mR^n))}^{4q}
\(\mE|\Lambda_{t}^{\e}-\Lambda_{t}^{0}|^{8q}\)^{\frac{1}{2}}
\(\mE(\Lambda_{T}^{\e})^{-2}\)^{\frac{1}{2}}\\
&\leq& C\(\mE|\Lambda_{t}^{\e}-\Lambda_{t}^{0}|^{8q}\)^{\frac{1}{2}}.
\de
 Next, we observe $|\Lambda_{t}^{\e}-\Lambda_{t}^{0}|$. By definitions of $\Lambda_{t}^{0}$ and $\Lambda_{t}^{\e}$, it holds that
\be
&&|\Lambda_{t}^{\e}-\Lambda_{t}^{0}|\no\\
&=&\Big|\exp\Big\{\int_{0}^{t}h^{i}(X_{s}^{\e},\sL^{\mP}_{X_{s}^{\e}})\dif V_{s}^{i} +\frac{1}{2}\int_{0}^{t}|h(X_{s}^{\e},\sL^{\mP}_{X_{s}^{\e}})|^{2}\dif s\Big\}-
\exp\Big\{\int_{0}^{t}h^{i}(\bar{X}_{s},\sL^{\mP}_{\bar{X}_{s}})\dif V_{s}^{i}\no\\ &&\quad+\int_{0}^{t}h^{i}(\bar{X}_{s},\sL^{\mP}_{\bar{X}_{s}})h^{i}(X_{s}^{\e},\sL^{\mP}_{X_{s}^{\e}})\dif s
-\frac{1}{2}\int_{0}^{t}|h(\bar{X}_{s},\sL^{\mP}_{\bar{X}_{s}})|^{2}\dif s\Big\}\Big|\no\\
&=&\Lambda_{t}^{0}\cdot\Big|\exp\Big\{\int_{0}^{t}\left(h^{i}(X_{s}^{\e},\sL^{\mP}_{X_{s}^{\e}})
-h^{i}(\bar{X}_{s},\sL^{\mP}_{\bar{X}_{s}})\right)\dif V_{s}^{i} \no\\ &&\quad+\frac{1}{2}\int_{0}^{t}\left(|h(\bar{X}_{s},\sL^{\mP}_{\bar{X}_{s}})|^{2}
-|h(X_{s}^{\e},\sL^{\mP}_{X_{s}^{\e}})|^{2}\right)\dif s\no\\
&&\quad+\int_{0}^{t}\left(|h(X_{s}^{\e},\sL^{\mP}_{X_{s}^{\e}})|^{2}-h^{i}(\bar{X}_{s},\sL^{\mP}_{\bar{X}_{s}})
h^{i}(X_{s}^{\e},\sL^{\mP}_{X_{s}^{\e}})\right)\dif s\Big\}-1\Big|.
\label{14}
\ee
Then, we deal with the integral $\int_{0}^{t}\left(h^{i}(X_{s}^{\e},\sL^{\mP}_{X_{s}^{\e}})-h^{i}(\bar{X}_{s},\sL^{\mP}_{\bar{X}_{s}})\right)\dif V_{s}^{i}$. From the isometric formula and $(\mathbf{H}_{h})$, it follows that
\ce
&&\mE\Big|\int_{0}^{t}\(h^{i}(X_{s}^{\e},\sL^{\mP}_{X_{s}^{\e}})
-h^{i}(\bar{X}_{s},\sL^{\mP}_{\bar{X}_{s}})\)\dif V_{s}^{i}\Big|^{2}\\
&=&\mE\Big[\int_{0}^{t}|h^{i}(X_{s}^{\e},\sL^{\mP}_{X_{s}^{\e}})
-h^{i}(\bar{X}_{s},\sL^{\mP}_{\bar{X}_{s}})|^{2}\dif s\Big]\\
&\leq&\mE\Big[\int_{0}^{t}L_{h}\Big(|X_{s}^{\e}-\bar{X}_{s}|^{2}
+\mW_2^{2}(\sL^{\mP}_{X_{s}^{\e}},\sL^{\mP}_{\bar{X}_{s}})\Big)\dif s\Big]\\
&\leq& C\mE\(\sup_{0\leq t\leq T}|X_{t}^{\e}-\bar{X}_{t}|^{2}\),
\de
which implies that
\ce
\lim_{\e\rightarrow0}\mE\Big|\int_{0}^{t}\(h^{i}(X_{s}^{\e},\sL^{\mP}_{X_{s}^{\e}})
-h^{i}(\bar{X}_{s},\sL^{\mP}_{\bar{X}_{s}})\)\dif V_{s}^{i}\Big|^{2}=0,
\de
and then
$$
\lim_{\e\rightarrow0}\int_{0}^{t}\left(h^{i}(X_{s}^{\e},\sL^{\mP}_{X_{s}^{\e}})
-h^{i}(\bar{X}_{s},\sL^{\mP}_{\bar{X}_{s}})\right)\dif V_{s}^{i}=0, \quad a.s..
$$
For the integral $\int_{0}^{t}(|h(\bar{X}_{s},\sL^{\mP}_{\bar{X}_{s}})|^{2}-|h(X_{s}^{\e},\sL^{\mP}_{X_{s}^{\e}})|^{2})\dif s$, by $(\mathbf{H}_{h})$ and the dominated convergence theorem, we obtain that
$$
\lim_{\e\rightarrow0}\int_{0}^{t}\left(|h(\bar{X}_{s},\sL^{\mP}_{\bar{X}_{s}})|^{2}-
|h(X_{s}^{\e},\sL^{\mP}_{X_{s}^{\e}})|^{2}\right)\dif s=0, \quad a.s..
$$
By the similar deduction to the above equality, one could get
$$
\lim_{\e\rightarrow0}\int_{0}^{t}\left(|h(X_{s}^{\e},\sL^{\mP}_{X_{s}^{\e}})|^{2}-
h^{i}(\bar{X}_{s},\sL^{\mP}_{\bar{X}_{s}})h^{i}(X_{s}^{\e},\sL^{\mP}_{X_{s}^{\e}})\right)
\dif s=0, \quad a.s..
$$
Thus, by taking the limit on both sides of (\ref{14}), it holds that
$$
\lim_{\e\rightarrow0}|\Lambda_{t}^{\e}-\Lambda_{t}^{0}|=0, \quad a.s..
$$

Also note that
\ce
|\Lambda_{t}^{\e}|^{8q}&=&\exp\Big\{8q\int_{0}^{t}h^{i}(X_{s}^{\e},\sL^{\mP}_{X_{s}^{\e}})\dif V_{s}^{i} +4q\int_{0}^{t}|h(X_{s}^{\e},\sL^{\mP}_{X_{s}^{\e}})|^{2}\dif s\Big\}\\
&=& \exp\Big\{\int_{0}^{t}8q h^{i}(X_{s}^{\e},\sL^{\mP}_{X_{s}^{\e}})\dif V_{s}^{i} -\frac{1}{2}\int_{0}^{t}|8q h(X_{s}^{\e},\sL^{\mP}_{X_{s}^{\e}})|^{2}\dif s\Big\}\\
&&\quad\cdot \exp\Big\{(32 q^2+4q)\int_{0}^{t}|h(X_{s}^{\e},\sL^{\mP}_{X_{s}^{\e}})|^{2}\dif s\Big\}\\
&\leq& \exp\Big\{\int_{0}^{t}8q h^{i}(X_{s}^{\e},\sL^{\mP}_{X_{s}^{\e}})\dif V_{s}^{i} -\frac{1}{2}\int_{0}^{t}|8q h(X_{s}^{\e},\sL^{\mP}_{X_{s}^{\e}})|^{2}\dif s\Big\}\exp\{CT\}.
\de
$$
|\Lambda_{t}^{0}|^{8q}\leq \exp\Big\{\int_{0}^{t}8q h^{i}(\bar{X}_{s},\sL^{\mP}_{\bar{X}_{s}})\dif V_{s}^{i} -\frac{1}{2}\int_{0}^{t}|8q h(\bar{X}_{s},\sL^{\mP}_{\bar{X}_{s}})|^{2}\dif s\Big\}\exp\{CT\}.
$$
Thus, by the dominated convergence theorem it holds that
\be
\lim_{\e\rightarrow0}A_{2}=0.
\label{15}
\ee

Finally, taking the limit on both sides of (\ref{11}) and inserting (\ref{13}) and (\ref{15}) in (\ref{11}), we have (\ref{rho}). The proof is complete.

\end{document}